\documentclass[11pt]{article}

\usepackage[margin=1in]{geometry}
\usepackage{amsmath}
\usepackage{amsthm}
\usepackage{bbm}
\usepackage{comment}
\usepackage{amssymb}
\usepackage{dsfont}
\usepackage{enumerate}
\usepackage{tikz}
\usepackage{hyperref}
\usepackage{sidecap}
\usepackage{caption,subcaption}
\usepackage[capitalise]{cleveref} 
\usepackage{setspace}

\newtheorem{theorem}{Theorem}
\newtheorem{proposition}[theorem]{Proposition}
\newtheorem{lemma}[theorem]{Lemma}
\newtheorem{corollary}[theorem]{Corollary}

\theoremstyle{definition}
\newtheorem{definition}[theorem]{Definition}
\numberwithin{theorem}{section} 

\renewcommand{\subset}{\subseteq}
\renewcommand{\tilde}{\widetilde}
\renewcommand{\hat}{\widehat}
\renewcommand{\epsilon}{\varepsilon}

\def\sign{{\rm sign}}
\def\sinc{{\rm sinc}}
\def\supp{{\rm supp}}

\def\wherespace{\quad\text{where}\quad}
\def\foreachspace{\quad\text{for each}\quad}
\def\forallspace{\quad\text{for all}\quad}
\def\andspace{\quad\text{and}\quad}

\def\<{\langle}
\def\>{\rangle}
\def\({\left(}
\def\){\right)}

\def\calC{\mathcal{C}}

\def\calF{\mathcal{F}}
\def\calG{\mathcal{G}}
\def\calH{\mathcal{H}}
\def\calI{\mathcal{I}}
\def\calJ{\mathcal{J}}
 
\def\calN{\mathcal{N}} 
\def\calP{\mathcal{P}}

\def\calU{\mathcal{U}}
\def\calV{\mathcal{V}}
\def\calW{\mathcal{W}}
\def\calX{\mathcal{X}}

\def\C{\mathbb{C}}

\def\N{\mathbb{N}}

\def\R{\mathbb{R}}
\def\S{\mathbb{S}}
\def\T{\mathbb{T}}
\def\Z{\mathbb{Z}}

\numberwithin{equation}{section}

\title{Nonharmonic multivariate Fourier transforms and matrices: condition numbers and hyperplane geometry}
\author{Weilin Li\footnote{CUNY City College. Email: wli6@ccny.cuny.edu}}

\begin{document}
\maketitle

\begin{abstract}
	Consider an operator that takes the Fourier transform of a discrete measure supported in $\mathcal{X}\subset[-\frac 12,\frac 12)^d$ and restricts it to a compact $\Omega\subset\R^d$. We provide lower bounds for its smallest singular value when $\Omega$ is either a closed ball of radius $m$ or closed cube of side length $2m$, and under different types of geometric assumptions on $\mathcal{X}$. We first show that if distances between points in $\mathcal{X}$ are lower bounded by a $\delta$ that is allowed to be arbitrarily small, then the smallest singular value is at least $Cm^{d/2} (m\delta)^{\lambda-1}$, where $\lambda$ is the maximum number of elements in $\mathcal{X}$ contained within any ball or cube of an explicitly given radius. This estimate communicates a localization effect of the Fourier transform. While it is sharp, the smallest singular value behaves better than expected for many $\mathcal{X}$, including when we dilate a generic set by parameter $\delta$. We next show that if there is a $\eta$ such that, for each $x\in\mathcal{X}$, the set $\mathcal{X}\setminus\{x\}$ locally consists of at most $r$ hyperplanes whose distances to $x$ are at least $\eta$, then the smallest singular value is at least $C m^{d/2} (m\eta)^r$. For dilations of a generic set by $\delta$, the lower bound becomes $C m^{d/2} (m\delta)^{\lceil (\lambda-1)/d\rceil }$. The appearance of a $1/d$ factor in the exponent indicates that compared to worst case scenarios, the condition number of nonharmonic Fourier transforms is better than expected for typical sets and improve with higher dimensionality. 
\end{abstract}

\medskip
\noindent
{\bf 2020 Math Subject Classification:} 15A12, 15A60, 42A05, 42A15, 65F22


\section{Introduction}

\label{sec:intro}

\subsection{Motivation}

Nonharmonic multidimensional Fourier transforms and matrices naturally appear in imaging applications, numerical schemes, and pure harmonic analysis. It is important to theoretically understand their condition number, since it provides insight into stability of inversion, properties of data processing algorithms, fundamental limits of recovery, stability of multivariate trigonometric interpolation, and more. Despite how frequently Fourier transforms appear and their connection to multivariate polynomial interpolation, there are many open questions about their extreme singular values. 

Let us first explain the above terminology and our motivation for our main questions. For a finite $\calX:=\{x_k\}_{k=1}^s\subset [-\frac 12,\frac 12)^d$, consider an arbitrary complex valued discrete measure $\mu$ supported in $\calX$, which can be written as
$$
\mu = \sum_{k=1}^s u_k \delta_{x_k}, \wherespace u\in \C^s. 
$$
Here, $\calX$ and $u$ are arbitrary. A vast collection of questions center around the recovery of $\mu$ when we are given Fourier data of $\mu$ on a compact $\Omega\subset\R^d$. The compactness assumption is essential to the formulation and applicability of this question because only a certain amount of information can be collected in practice due to physical limitations and the complexity of computational schemes scale in the number of samples. This problem is interesting when $\calX$ is arbitrary and can be nonuniform with geometric properties. 

More specifically, we first consider a continuous model. Suppose we take the Euclidean Fourier transform of $\mu$ and restrict it to $\Omega$, which we denote by $\hat\mu \, |_\Omega$. We are lead to a {\it nonharmonic Fourier transform} $\calF_{\Omega,\calX}\colon \C^s\to L^2(\Omega)$ such that 
\begin{equation}
	\label{eq:Fomegadef}
	\calF_{\Omega,\calX}u:= \hat\mu \, \big|_{\Omega}\andspace
	\hat\mu(\omega) := \int_{[-\frac 12,\frac 12)^d} e^{-2\pi i \omega \cdot x} \, d\mu(x) = \sum_{k=1}^s u_k e^{-2\pi i\omega\cdot x_k}.  
\end{equation}
We also consider a discrete model whereby Fourier data is sampled on a finite set. Since we work with the normalization $\calX\subset[-\frac 12,\frac 12)^d$, the canonical discrete version corresponds to sampling on $\Omega \cap\Z^d$. In this case, we need to interpret $\calX$ as a subset of the torus $\T^d:=(\R\setminus\Z)^d$ and $\hat\mu$ will be the Fourier transform on $\T^d$. We are led to a {\it nonharmonic multivariate Fourier matrix}, 
\begin{equation}
	\Phi_{\Omega,\calX}:=\Big[ e^{-2\pi i \omega\cdot x} \Big]_{\omega\in \Omega\cap\Z^d, \, x\in \calX} \andspace \Phi_{\Omega,\calX}u = \Big[\hat\mu (\omega) \Big]_{\omega\in\Omega\cap\Z^d}.
\end{equation}

In these definitions, we place no assumptions on $\Omega$, since there is a rich variety of sampling sets that are pertinent in applications. For theoretical analysis, we will consider the most natural cases such as when $\Omega$ is a ball or cube. They are both nonharmonic Fourier transforms because the functions $\{\omega\mapsto e^{2\pi i \omega\cdot x_k}\}_{k=1}^s$ are not $L^2$ orthogonal, except for special cases, such as when $\calX$ consists of uniformly spaced points. 

The focus of this paper is on the minimum singular values of these operators. They are defined in the usual way,
$$
\sigma_{\min}(\calF_{\Omega,\calX}):= \inf_{u\in \C^s, \, |u|_2=1} \|\calF_{\Omega,\calX} u\|_{L^2(\Omega)} \andspace
\sigma_{\min}(\Phi_{\Omega,\calX}):= \min_{u\in \C^s, \, |u|_2=1} |\Phi_{\Omega,\calX} u|_2. 
$$
In the univariate case, it is well understood that the smallest singular value greatly depends on $m$, $s$, and configuration of $\calX$, see \cite{gautschi1963inverses,bazan2000conditioning,moitra2015matrixpencil,aubel2019vandermonde,li2021stable,batenkov2020conditioning,kunis2020smallest,batenkov2021spectral,li2023multiscale}. On the other hand, the maximum singular value is not as interesting to study. Here and throughout this paper, we let $|\Omega|$ be the Lebesgue measure of $\Omega$, while we let $|\Omega|_*=|\Omega\cap\Z^d|$ be the number of multi-integers contained in $\Omega$. It is not hard to show that 
\begin{equation}
	\label{eq:sigmax}
	\sqrt{|\Omega|} \leq \sigma_{\max}(\calF_{\Omega,\calX})\leq \sqrt{s|\Omega|} \andspace
	\sqrt{|\Omega|_*} \leq \sigma_{\max}(\Phi_{\Omega,\calX})\leq \sqrt{s|\Omega|_*}.
\end{equation}
From this discussion, we see that any lower bound for the smallest singular value will provide an upper bound for the condition number. 

We study both Fourier operators rather than just one for several reasons. The continuous operator $\calF_{\Omega,\calX}$ is more convenient to work with, since it enjoys additional properties (e.g., certain symmetry and dilation relationships) that $\Phi_{\Omega,\calX}$ does not, due to the latter having discrete set of frequencies. The former is attractive from a theoretical perspective, while the latter is more relevant in applications. By studying both operators, we can disentangle the effects of band-limitation that is present in both $\calF_{\Omega,\calX}$ and $\Phi_{\Omega,\calX}$, versus the combined effects of band-limitation and sampling which is only present in $\Phi_{\Omega,\calX}$. Finally, we might not necessarily want to discretize $\Omega$ via $\Z^d$ and would prefer other sampling methods (e.g., random or oversampling), so it is convenient to analyze the properties of $\calF_{\Omega,\calX}$ directly.

\subsection{Some geometric concepts} 

The main techniques and ideas used in this paper depend on the choice of metric and different geometric descriptions. As remarked earlier, we will generally analyze the case when $\Omega$ is a ball or cube. For $p\in [1,\infty]$, the closed $\ell^p$ ball in $\R^d$ with radius $m$ centered at the origin is denoted
$$
\Omega_m^p:=\{x\in \R^d\colon |x|_p\leq m\}.
$$
We can interpret $2m$ as the bandwidth of our measurements. Since we will pay special attention to the ball and cube cases, in order to simplify notation when appropriate, we let 
$$
B_m:=\Omega_m^2
\andspace Q_m:=\Omega_m^\infty. 
$$

When working with $\calF_{\Omega,\calX}$, we use the usual $\ell^p(\R^d)$ norm $|\cdot|_p$ to quantify distances between pairs in $\calX\subset\R^d$. When dealing with $\Phi_{\Omega,\calX}$, we use the $\ell^p(\T^d)$ norm instead, since $\calX\subset\T^d$ is treated as a periodic set, and we use the same notation $|\cdot|_p$. In either case, we define the $\ell^p$ minimum separation of $\calX$ as 
\begin{equation}
	\Delta_p(\calX) = \min_{x,x'\in \calX} |x-x'|_{p}. 
\end{equation}
From duality, if the Fourier transform is restricted to a $\ell^p$ ball, it is sometimes natural to use the $\ell^{p'}$ metric on $\calX$, where $p'$ is the H\"older dual of $p$. 

We define several quantities that depend on the behavior of $\calX$ locally. For a positive $\tau$, we denote the $\ell^p$ {\it neighborhood set} of a $x\in \calX$ as
$$
\calN_p(x,\tau,\calX)
:=\{x'\in \calX\colon |x'-x|_p\leq \tau\}. 
$$
In other words, $\calN_p(x,\tau,\calX)$ is the intersection of $\calX$ with a closed $\ell^p$ ball of radius $\tau$ centered at $x$. We use the following definition in \cite{li2023multiscale}.
\begin{definition}
	The $\tau$ {\it local sparsity} of $\calX$ is defined as
	\begin{equation}
		\label{eq:deflocalsparsity}
		\nu_p(\tau,\calX)
		:= \max_{x\in \calX} \, |\calN_p(x,\tau,\calX)|
		\andspace \nu_p(\tau,\emptyset):=0. 
	\end{equation}
\end{definition}

In other words, the local sparsity $\nu_p(\tau,\calX)$ is the maximum number of $x'\in\calX$ contained in any neighborhood set at scale $\tau$. Many of our results require some condition on the local sparsity versus $m$. To give more concrete examples of sets that will satisfy some of our main theorems, we introduce multidimensional versions of clumps.

\begin{definition}
	We say $\calX$ {\it consists of clumps} with parameters $(p,\tau,\lambda)$ if the following hold. There exist $r\in \N_+$ and disjoint subsets $\calC_1,\dots, \calC_r$ of $\calX$ called {\it clumps} such that $\calX = \calC_1\cup \cdots \cup \calC_r$, the $\ell^p$ distance between any two elements in the same clump is bounded above by $\tau$, and {\color{blue} $\max_{1\leq k\leq r}$}  $|\calC_k|=\lambda$. Further, if $r>1$, assume that the $\ell^p$ distance between different $\calC_j$ and $\calC_k$ exceeds $\tau$.  
\end{definition}

Aside from clumps, a more complicated geometric description can be given in terms of hyperplane decompositions.  We say a set $\calX$ {\it consists of $r$ hyperplanes relative to a $x\in \calX$} (called the {\it reference point}) if there exist hyperplanes $\calH_1,\dots, \calH_r\subset\R^d$ that do not intersect $x$ and
\begin{equation}
	\label{eq:hyperdecomp}
	\calX\setminus \{x\}=\calH_1\cup \cdots\cup\calH_r.
\end{equation}
The same definition appeared in Chung-Yao \cite{chung1977lattices}, which they called the geometric characterization. This decomposition is generally not unique since there may be many valid choices of hyperplanes. It is important to remark that this definition allows for $\calX \cap \calH_j$ and $\calX \cap\calH_k$ have non-empty intersection for $j\not=k$. 

\begin{figure}[ht]
	\centering
	\begin{tikzpicture}[scale=.75]
		\filldraw[color=black!10] (0,0) circle (122pt);
		\draw[thick] (0,0) circle (122pt);
		\draw[thick, dashed] (-3,3)--(4,1.25);
		\draw[thick, dashed] (-2,-3.75)--(-2,3.75);
		\filldraw (0,0) circle (2pt) node[anchor=north]{$x_1$};
		\filldraw (-1,2.5) circle (2pt) node[anchor=south]{$x_2$};
		\filldraw (3,1.5) circle (2pt) node[anchor=south]{$x_3$};
		\filldraw (-2,-2.5) circle (2pt) node[anchor=east]{$x_4$};
		\draw[thick,->] (0,0)--(3,-3);
		\draw (2,-1) node[anchor=north]{$\tau$};
		\filldraw (16,0) circle (2pt) node[anchor=north]{$x_6$};
		\filldraw (13,-1) circle (2pt) node[anchor=north]{$x_7$};
		\draw[thick,->] (0,0)--(-2,0) node[anchor=east]{$\calH_2$};
		\draw (-1,0) node[anchor=south]{$\eta_2$};
		\draw[thick,->] (0,0)--(0.52,2.08) node[anchor=south]{$\calH_1$};
		\draw (0.2,1) node[anchor=west]{$\eta_1$};
		\filldraw (8,-2) circle (2pt) node[anchor=north]{$x_5$};
		\draw[thick,->]  (0,0)--(8,-2);
		\draw (5.5,-1.25) node[anchor=south]{$>\tau$};
	\end{tikzpicture}
	\caption{Example of local hyperplane decomposition with reference point $x_1$. The closed ball of radius $\tau$ centered at $x_1$ is shown in gray, and it contains $x_1$, $x_2$, $x_3$, and $x_4$, which is the neighborhood set $\calN_2(x_1,\tau,\calX)$. The hyperplane $\calH_1$ contains $x_2$ and $x_3$, while the hyperplane $\calH_2$ contains just $x_4$. The distance between these hyperplanes to $x_1$ are $\eta_1$ and $\eta_2$, respectively.}
	\label{fig:localhyperplane}
\end{figure}
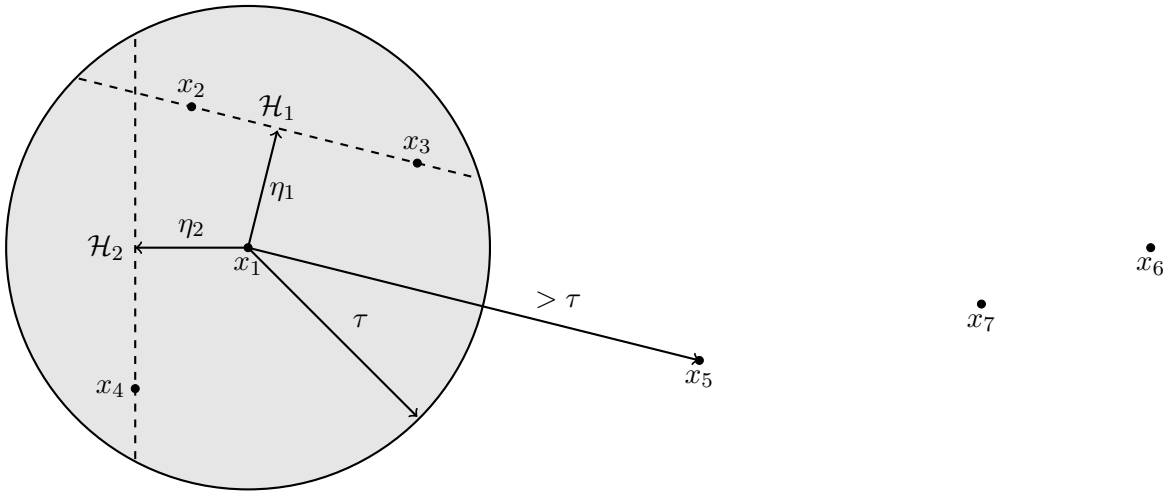

The above hyperplane decomposition is global by definition. We need to localize it so that only distances less than a prescribed scale influence the number of hyperplanes required. See Figure \ref{fig:localhyperplane} for an illustration of a local hyperplane decomposition. 

\begin{definition}
	We say $\calX$ {\it locally consists of hyperplanes} (or {\it admits a local hyperplane decomposition}) with parameters $(p,\tau,r,\eta)$ if for each reference point $x\in \calX$, the neighborhood set $\calN_p(x,\tau,\calX)$ consists of at most $r$ hyperplanes relative to $x$, and the $\ell^2$ distance between $x$ and each hyperplane is lower bounded by $\eta$. 
\end{definition}

\subsection{Summary of main results and organization}

The main contributions this paper consist of accessible and accurate estimates for $\sigma_{\min}(\calF_{\Omega_m^p,\calX})$ and $\sigma_{\min}(\Phi_{\Omega_m^p,\calX})$ with $p\in \{2,\infty\}$, in terms of reasonably simple geometric descriptions of $\calX$. We present estimates for three drastically different settings. In this summary section, we state everything in an informal manner and make over simplifications in order to present the general flavor of results. All constants that appear in this section do not depend on $m$ or the minimum separation of $\calX$, but may depend on other quantities such as $p$, cardinality of $\calX$, local sparsity of $\calX$ at some scale, and dimension $d$.

In \cref{sec:wellsep}, we examine the {\it well-separated case}, which loosely corresponds to the setting where $\Delta_p(\calX)$ is sufficiently large compared to $\frac 1m$. \cref{thm:wellsepball,thm:wellsepcube} provide an explicit $C,c>0$ that do not depend on $m$ such that if $\Delta_p(\calX)\geq \frac Cm$, then  
$$
\sigma_{\min}(\calF_{\Omega_m^p,\calX}) \geq c\sqrt{|\Omega_m^p|} \andspace
\sigma_{\min}(\Phi_{\Omega_m^p,\calX}) \geq c \sqrt{|\Omega_m^p|_*}. 
$$
Consequently, in the well-separated case, both operators have condition numbers that are upper bounded by constants that do not depend on $m$. These results are proved by using a well-known technique that exploits properties of special functions. We will provide some comparisons with prior work in that section.

There are many applications where a strong separation condition is violated. In \cref{sec:sr1}, we focus on the {\it super-resolution regime}, which allows for $\Delta_{p'}(\calX)$ to be arbitrarily small. As the separation goes to zero, the smallest singular values of both Fourier operators tend to zero, and our primary goal is to derive non-asymptotic and sharp lower bounds. \cref{thm:srball,thm:srcube} are multiscale estimates which show that the distances between elements in $\calX$ at small scales predominately influence these operators' smallest singular values. We specialize these theorems to clumps in order to give more interpretable results. \cref{cor:clumpsball,cor:clumpscube} show that if $\calX$ consists of clumps with parameters $(p,\tau,\lambda)$ where $\tau \gtrsim \frac \lambda m$ and $\Delta_{p'}(\calX)$ is sufficiently small, then
$$
\sigma_{\min}(\calF_{\Omega_m^p,\calX}) 
\gtrsim \sqrt{|\Omega_m^p|} \, \big(m \Delta_{p'}(\calX)\big)^{\lambda-1} 
\andspace
\sigma_{\min}(\Phi_{\Omega_m^p,\calX}) 
\gtrsim \sqrt{|\Omega_m^p|_*} \, \big(m \Delta_{p'}(\calX) \big)^{\lambda-1}. 
$$
These estimates are sharp in $m$, $\lambda$, and $\Delta_p(\calX)$. 

In \cref{sec:sr2}, we continue with the super-resolution regime, but we bound their singular values in terms of a hyperplane decomposition of $\calX$. This material is motivated by numerical experiments that indicate the previously stated inequalities are not sharp if we dilate generic $\calX$ by $\delta$. To partially explain this behavior, \cref{thm:hyperball,thm:hypercube} show that if $\calX$ locally consists of hyperplanes with parameters $(p,\frac {\nu_p(\tau,\calX)} m,r,\eta)$ for an appropriate $\tau$, then
$$
\sigma_{\min}(\calF_{\Omega_m^p,\calX}) 
\gtrsim \sqrt{|\Omega_m^p|} \, (m\eta)^{r}. 
$$
When specializing this to the case where $\calX$ consists of at most $\lambda$ generic points that are close together, then $r=\lceil \frac{\lambda-1}{d} \rceil$. This exponent is sharp for small enough $\lambda$ depending on $d$; for instance, when $d=2$, it is sharp for $\lambda \in \{1,\dots,5,7\}$. 

By comparing the above findings for arbitrary versus generic sets, we see that the condition number behaves better than expected for sets with favorable local hyperplane structure. For dilations of a generic set, the appearance of a $1/d$ factor in the expression for $r$ also indicates that there is a blessing of dimensionality. Consequently, inversion of nonharmonic Fourier transforms is more stable than what one initially suspects from just worst case analysis and results for the one-dimensional problem.

\subsection{Preliminary notation}

The torus is defined as $\T:=\R/\Z$ and is equipped with the metric $d_{\T}(s,t):=|s-t|_{\T}:= \min_{n\in \Z} |s-t+n|$. We normally identify $\T^d$ with $[-\frac 1 2, \frac 1 2)^d$. We use $|\cdot|_p$ for the $\ell^p$ norm on $\R^d$ or $\T^d$. For a Lebesgue measurable $A\subset \R^d$ or $A\subset\T^d$, let $L^p(A)$ be the set of complex-valued functions $f$ defined on that are Lebesgue measurable such that $|f|^p$ is integrable. We use $\|\cdot\|_{L^p}$ for the $L^p$ norm and let $p'$ be H\"older conjugate of $p$, namely $\frac 1 p + \frac 1{p'}=1$ with the usual convention that $0=\frac 1 \infty$. We use the same notation for the Fourier transform on $\R^d$ and $\T^d$. For a suitable $f$, its Fourier transform is denoted $\hat f$, where $\hat f(\xi):=\int f(x) e^{-2\pi i \xi\cdot x}\, dx$. Let $f|_A$ denote the restriction of a function to $A$ and $\mathbbm{1}_A$ be the indicator function of $A$. The additive sum of two sets $A,B$ is denoted $A+B$. We use the notation $x\lesssim_{a,b} y$ to mean that there is a constant $C>0$ that only depends on $a$ and $b$ such that $x\leq Cy$. We also use the notation $x\asymp_{a,b} y$ to mean that $x\lesssim_{a,b} y$ and $y\lesssim_{a,b} x$ simultaneously hold. Additional definitions and notation will be introduced in later sections.

\section{Well-separated case} 
\label{sec:wellsep}

In this section, we state our main results for the well-separated case and introduce all the necessary proof ingredients. The primarily technical tool that we will use are properties of certain functions. We begin with an abstract lemma, which is proved in \cref{proof:poisson}. 

\begin{lemma}
	\label{lem:poisson}
	Let $\Omega\subset\R^d$ be compact, $p\in [1,\infty]$, and $\delta \in (0,\frac 1 2)$. Suppose there exist $\varphi$ and $\psi$ that are continuous on $\R^d$, supported in $\Omega_\delta^p$, $\hat\varphi, \hat\psi\in L^1(\R^d)$, and $\hat \psi(\xi)\leq \mathbbm{1}_\Omega(\xi)\leq \hat\varphi(\xi)$ for all $\xi\in \R^d$. For any set $\calX\subset [-\frac 12,\frac 12)^d$ such that $\Delta_p(\calX)>\delta$, we have 
	\begin{align*}
		\sqrt{\psi(0)}
		\leq \sigma_{\min}(\calF_{\Omega,\calX})
		&\leq \sigma_{\max}(\calF_{\Omega,\calX})
		\leq \sqrt{\varphi(0)}, \\
		\sqrt{\psi(0)}
		\leq \sigma_{\min}(\Phi_{\Omega,\calX})
		&\leq \sigma_{\max}(\Phi_{\Omega,\calX})
		\leq \sqrt{\varphi(0)}.
	\end{align*}
\end{lemma}

Any $\varphi$ and $\psi$ that satisfy the properties of this lemma are loosely referred to as a {\it majorant} and {\it minorant} of $\mathbbm{1}_\Omega$, respectively. The main difficulties are usually with the existence of a minorant, and we say $\psi$ is {\it nontrivial} if $\psi(0)>0$. A nontrivial minorant does not always exist. For instance, the uncertainty principle for the Fourier transform tells us that a function cannot be too concentrated in both time and frequency. It follows from this principle that $\psi$ will never exist if $\delta$ is too small depending on $\Omega$. 

Beurling-Selberg constructed a minorant $\psi$ for an interval $\Omega=[-m,m]$, see \cite{selberg1989collected} for details and history. This result has found many applications in number theory \cite{montgomery1978analytic} and one-dimensional super-resolution \cite{moitra2015matrixpencil,liao2016music,aubel2019vandermonde,li2020super}. There is a complete answer for when $\Omega$ is the unit ball \cite{holt1996beurling,gonccalves2018note}, while there are partial results for the cube \cite{barton1999analogs,carruth2019extremal,carruth2022beurling}. To our best knowledge, this question is completely open for other $\ell^p$ balls and other convex sets.

We first state results for the cube case since they are easier to describe. In this paper, we use a minorant in Barton's PhD thesis that was denoted $M_3^-$.

\begin{theorem}[Barton, \cite{barton1999analogs}, Section 2.5]
	\label{thm:barton}
	For any $\delta>0$ and integer $d\geq 2$, the following hold. There exist $\psi$ and $\varphi$ that are uniformly continuous on $\R^d$, supported in $Q_\delta$, $\hat\varphi,\hat\psi\in L^1(\R^d)$, $\hat\psi(\xi)\leq \mathbbm{1}_{Q_1}(\xi)\leq \hat\varphi(\xi)$ for all $\xi\in\R^d$, and 
	\begin{align*}
		\int_{\R^d} (\hat\varphi-\mathbbm{1}_{Q_1})
		=\int_{\R^d} (\mathbbm{1}_{Q_1}-\hat\psi)
		&=\( \big(1+\tfrac 1 {2\delta }\big)^d-1 \) \, |Q_1|.
	\end{align*}
\end{theorem}

While $\delta$ is arbitrary in the statement of this result, there is an implicit restriction. If $\delta$ is too small, then we cannot guarantee that $\psi(0)>0$. However, if $\delta \geq \frac{d}{2\log 2}$, then using the inequality $1+t<e^t$ for all $t>0$, we see that 
\begin{align}
	\label{eq:psi0barton}
	\psi(0)
	&=|Q_1|-\int_{\R^d} (\mathbbm{1}_{Q_1}-\hat\psi) 
	= \( 2-\big(1+\tfrac 1 {2\delta }\big)^d\) \, |Q_1|
	> ( 2- e^{d/(2\delta)}) |Q_1| 
	\geq 0, \\
	\label{eq:psi0barton2}
	\varphi(0)
	&=|Q_1|+\int_{\R^d} (\hat\varphi-\mathbbm{1}_{Q_1}) 
	= \big(1+\tfrac 1 {2\delta }\big)^d\, |Q_1|
	< e^{d/(2\delta)} |Q_1|.   
\end{align}
By combining \cref{lem:poisson} and \cref{thm:barton}, we obtain the following theorem for the cube, which is proved in \cref{proof:wellsepcube}. 

\begin{theorem}[Samples in a cube and well separated nodes]
	\label{thm:wellsepcube}
	For any integer $d\geq 2$, $m\geq 1$, $\beta \geq \frac 1 {2\log 2}$, and $\calX\subset [-\frac 12,\frac 12)^d$ such that $\Delta_\infty(\calX) \geq \frac {\beta d} m$, we have
	\begin{align*}
		\sqrt{(2- e^{1/(2\beta)} )|Q_m|}
		&\leq \sigma_{\min}(\calF_{Q_m,\calX})
		\leq \sigma_{\max}(\calF_{Q_m,\calX})
		\leq \sqrt{e^{1/(2\beta)}|Q_m|}, \\
		\sqrt{(2- e^{1/(2\beta)} )|Q_m|_*}
		&\leq \sigma_{\min}(\Phi_{Q_m,\calX})
		\leq \sigma_{\max}(\Phi_{Q_m,\calX})
		\leq \sqrt{e^{1/(2\beta)}|Q_{m}|_*}. 
	\end{align*}
\end{theorem}

Notice that the $\ell^\infty$ separation condition for $\calX$ is inversely proportional to $m$, which is expected from scaling arguments. Additionally, the constant $\beta d$ grows linearly in the dimension $d$, but it is not clear whether this growth in $d$ is necessary. 

Let us compare \cref{thm:wellsepcube} with Kunis-Nagel-Strotmann \cite[Theorem  4.2]{kunis2022multivariate}. They parameterized the frequency set by $\{0,\dots,N-1\}^d$, while we use $Q_m=\{-m,\dots,m\}^d$. The singular values of $\Phi_{\Omega,\calX}$ are invariant under multi-integer shifts of $\Omega$, so our settings are equivalent whenever we enforce the relationship $N=2m+1$. They proved that $\Delta_2(\calX)> \frac{4d}{2m+1}$ implies $\sigma_{\min}(\Phi_{Q_m,\calX})>0.9 \,\sqrt{N^d}=0.9\sqrt{|Q_{m}|_*}$. To get the same conclusion, we pick $\beta = (2\log(2-0.81))^{-1}\approx 0.7671$ in \cref{thm:wellsepcube} to see that $\Delta_2(\calX)\geq \frac{\beta d}{m}$ implies $\sigma_{\min}(\Phi_{Q_m,\calX}) \geq 0.9 \, \sqrt{|Q_m|_*}$. Both results require separation conditions that grow linearly in $d$, but our universal constant is smaller due to Barton.

Let us move onto the ball case. By exploiting rotational symmetry, Holt-Vaaler \cite{holt1996beurling} reduced the multidimensional problem to a single dimensional one. Building upon this result, Gon\c calves \cite{gonccalves2018note} derived an explicit expression for the smallest possible $\delta$ for which there exists a nontrivial minorant for the unit ball in arbitrary dimensions. 

Let $J_\nu$ be the Bessel function of the first kind of order $\nu$, as defined in Watson \cite{watson1922treatise}, and $j_{\nu,k}$ denote its $k$-th positive zero. Although Bessel functions can be defined on subsets of $\C$ and for complex $\nu$, we only need their restrictions to $\R$ and for real $\nu$. For $\alpha\in (\frac{j_{d/2-1,1}}\pi,\frac{j_{d/2,1}}\pi)$, let
\begin{equation}
	\label{eq:Calphadef}
	c(\alpha):=\frac{1}{|\S^{d-1}|} \(\frac 2 \alpha\)^d \frac{\gamma(\alpha)}{1+\gamma(\alpha)/d}
	\andspace \gamma(\alpha):= \frac{-\pi \alpha J_{d/2-1}(\pi \alpha)}{J_{d/2}(\pi \alpha)}. 
\end{equation}

\begin{theorem}[Gon\c calves \cite{gonccalves2018note}, Theorem 1]
	\label{thm:goncalves}
	Let $d\geq 2$ be an integer and $\alpha\in (\frac{j_{d/2-1,1}}\pi, \frac{j_{d/2,1}}\pi)$. Then there exists a uniformly continuous $\psi\colon\R^d\to\R$ supported in $B_\alpha$ such that $\hat\psi\in L^1(\R^d)$, $\hat\psi(\xi)\leq \mathbbm{1}_{B_1}(\xi)$ for all $\xi \in\R^d$ and 
	$
	\psi(0)= c(\alpha).
	$
	If $\alpha \leq \frac{j_{d/2-1,1}}\pi$, there is no $\psi$ satisfying these properties with $\psi(0)>0$.
\end{theorem}

Asymptotics for zeros of Bessel functions show that ${j_{d/2-1,1}}$ and ${j_{d/2,1}}$ grow linearly in $d$, see \cref{table:prelim} for a list of values for dimensions $d\in \{2,\dots,10\}$. For reasons that will be apparent momentarily, it will be convenient to extend the definition of $c(\alpha)$ so that 
\begin{equation}
	\label{eq:Calphadef2}
	c(\alpha) 
	:= \frac{1}{|\S^{d-1}|} \(\frac {2\pi } {j_{d/2,1}}\)^d \quad\text{if}\quad \alpha \geq \frac{j_{d/2,1}}\pi.
\end{equation}

\begin{table}[ht]
	\centering
	\begin{tabular}{|c|c|c|c|c|c|c|c|c|c|}
		\hline
		$d$ &2 &3 &4 &5 &6 &7 &8 &9 &10\\ \hline
		$\frac{j_{d/2-1,1}}\pi$ &0.7655 &1.0000 &1.2197 &1.4303 &1.6347 &1.8346 &2.0309 &2.2243  &2.4154 \\ \hline 
		$\frac{j_{d/2,1}}\pi$ &1.2197 &1.4303 &1.6347 &1.8346 &2.0309 &2.2243 &2.4154 &2.6046 &2.7920\\ \hline 
		$\sqrt{c(\frac{j_{d/2,1}}\pi)}$ &0.5220   &0.3947    &0.3033    &0.2357    &0.1848
		&0.1459    &0.1158    &0.0924    &0.0740 \\ \hline 
	\end{tabular}
	\caption{Numerically computed $\frac{j_{d/2-1,1}}\pi$, $\frac{j_{d/2,1}}\pi$, and $\sqrt{c(\frac{j_{d/2,1}}\pi)}$ for $d\in \{2,3,\dots,10\}$.}
	\label{table:prelim}
\end{table}

\begin{theorem}[Samples in a ball and well separated nodes]
	\label{thm:wellsepball}
	For any integer $d\geq 2$, $m>0$, $\alpha>\frac{j_{d/2-1,1}}\pi$, and $\calX\subset [-\frac 12,\frac 12)^d$ such that $\Delta_2(\calX)\geq \frac {\alpha} m$, we have 
	$$
	\sigma_{\min}(\calF_{B_m,\calX})
	\geq \sqrt{c(\alpha)m^d} \andspace \sigma_{\min}(\Phi_{B_m,\calX})
	\geq \sqrt{c(\alpha)m^d}. 
	$$
\end{theorem}

This theorem is proved in \cref{proof:wellsepball}. Just like in \cref{thm:wellsepcube}, the separation criterion is proportional to $\frac dm$. Unlike our result for the cube which yielded dimension independent constants, the quantity $c(\alpha)$ behaves poorly in high dimensions and one easily sees that for any choice of $\alpha$, we have $c(\alpha)\to 0$ super exponentially in $d$, see the last row of \cref{table:exponent}. To obtain improved constants, we can enforce a stronger condition on the minimum separation. If we had used the minorant constructed in \cite{holt1996beurling} instead of \cite{gonccalves2018note}, we could have shown that $\Delta_2(\calX)\geq \frac{cd}{m}$ for $c\approx 1.3356$ ensures that 
$\sigma_{\min}(\calF_{B_m,\calX}) \geq \frac 12 \sqrt{|B_m|}$ and likewise for $\sigma_{\min}(\Phi_{B_m,\calX})$. 

We have decided not to include this alternative result for several reasons. In low dimensions, it gives comparable estimates to \cref{thm:wellsepball} while requiring a considerably stronger separation requirement on $\calX$. In high dimensions, it greatly improves the constants, but we found that higher dimensional applications, such as nonuniform discrete Fourier transforms \cite{potts2001fast} and MIMO radar \cite{li2008mimo}, typically use samples in a cube rather than in a ball. Finally, the main error estimate in \cite{holt1996beurling} contains a weighted integral involving products of Bessel functions, so our proof of the alternative result requires a considerable detour into Bessel functions and explicit control over their asymptotic expansions. 

Although \cref{thm:wellsepball,thm:wellsepcube} are new, they are not the most novel aspects of this paper. They are essentially folklore since they result from combining the extremal function literature with a Poisson summation argument, such as \cref{lem:poisson}. This connection is well known to researchers working on the mathematics of super-resolution. The multivariate results presented in this section are generalizations of one-dimensional estimates that have appeared in \cite{moitra2015matrixpencil,liao2016music,aubel2014super} which show that a separation of $\frac C m$ implies the smallest singular value of a one-dimensional Fourier matrix (i.e., Vandermonde) is at least $c\sqrt m$. In comparison to prior work for the multidimensional case, we have made some improvements to explicit constants, which are important, since they will appear in a bootstrapping argument given in the next section for the super-resolution regime. 

There is a related but not directly comparable result. For $d=2$, Chen-Moitra \cite{chen2021algorithmic} studied how small $\Delta_2(\calX)$ can be so that $\sigma_{\min}(\calF_{B_{1/2},\calX})$ is not exponentially small in $s$. By selecting $\calX$ to be $s$ points on a hexagonal lattice of $\R^2$ and letting $s$ grow, they showed that $\Delta_2(\calX)\geq \sqrt{4/3}\approx 1.1547$ is required for stable recovery of $\mu$, in some appropriate sense defined in their paper. This construction can be adapted to derive an upper bound (as opposed to lower) for $\sigma_{\min}(\calF_{B_{1/2},\calX})$. Their conclusion that a minimum separation of $\sqrt{4/3}$ is necessary requires allowing $s$ to arbitrarily large. We will see in the next section that provided $s$ is fixed, the minimum separation can be arbitrarily small.

\section{Super-resolution regime and arbitrary nodes}

\label{sec:sr1}

In the super-resolution regime, a primary method that we will use to lower bound $\sigma_{\min}(\calF_{\Omega,\calX})$ and  $\sigma_{\min}(\Phi_{\Omega,\calX})$ is through a ``dual" relationship with minimum norm interpolation. This duality was introduced in \cite[Proposition 2.12]{li2021stable} for dealing with $\Phi_{\Omega,\calX}$ in the one dimensional case. It readily generalizes to higher dimensions and to $\calF_{\Omega,\calX}$ with minor modifications. This dual relationship was further combined with density ideas and sparsity decompositions in \cite{li2023multiscale}, which significantly improved existing estimates for one-dimensional Fourier matrices. 

We first collect these ideas before moving onto the main theorems. While we provide some exposition, it is rather brief; for those interested in the intuition and concepts behind these techniques, we refer the reader to \cite{li2023multiscale}, which explains them in greater detail and in the simpler one-dimensional case. Recall that throughout, we assume $\Omega\subset\R^d$ is compact. We let  $\calW(\Omega)$ denote the space of bandlimited functions or the
Paley-Wiener space; that is, a complex-valued $f$ belongs to $\calW(\Omega)$ if and only if there is a $F\in L^2(\Omega)$ such that 
$$
f(x)=\int_\Omega F(\omega) e^{2\pi i \omega\cdot x}\, d\omega \foreachspace x\in \R^d.
$$
Not surprisingly, $\calW(\Omega)$ is the natural space of functions associated to $\calF_{\Omega,\calX}$. When dealing with $\Phi_{\Omega,\calX}$, we examine $\calP(\Omega)$, the set of all trigonometric polynomials whose Fourier coefficients are supported in the finite set $\Omega \cap \Z^d$; that is, a complex-valued $f$ belongs to $\calP(\Omega)$ if and only if there is a $F\colon \Omega\cap \Z^d\to\C$ such that
$$
f(x)=\sum_{\omega \in \Omega\cap \Z^d} F(\omega) e^{2\pi i \omega \cdot x} \foreachspace x\in \T^d. 
$$
In either case, we identify $F(\omega)$ with $\hat f(\omega)$.

\begin{proposition}[Duality principle]
	\label{prop:duality}
	Suppose $\calX=\{x_k\}_{k=1}^s\subset [-\frac 12,\frac 12)^d$ and $\Omega\subset\R^d$ is compact. For any $\epsilon \in \C^s$ with $|\epsilon|_2<1$ and unit $\ell^2$ norm vector $v\in \C^s$, the following hold. If there exists $f\in \calW(\Omega)$ such that $f(x_k)=v_k+\epsilon_k$, then we have
	$$
	\|\calF_{\Omega,\calX} v\|_{L^2(\Omega)} 
	\geq (1-|\epsilon|_2) \, \|f\|_{L^2(\R^d)}^{-1}. 
	$$
	If there exists $f\in \calP(\Omega)$ such that $f(x_k)=v_k+\epsilon_k$, then we have
	$$
	|\Phi_{\Omega,\calX} v|_2 
	\geq (1-|\epsilon|_2) \, \|f\|_{L^2(\T^d)}^{-1}. 
	$$
\end{proposition}

This is a straightforward adaptation of \cite[Proposition 2.12]{li2021stable}, but we include its proof in \cref{proof:duality} for completeness. Typically, we will use this proposition for the case when $\epsilon=0$ and $v$ is a right singular vector of $\calF_{\Omega,\calX}$ (or $\Phi_{\Omega,\calX}$) corresponding to its smallest singular value. It is prudent to mention that this proposition does not make any claims about the relationship between singular values $\calF_{\Omega,\calX}$ and $\Phi_{\Omega,\calX}$ because they in general have different right singular vectors, so this statement would be applied to different $v$'s. 

We refer to this proposition as the {\it duality principle}, since it provides a connection between the singular values of these operators to minimum $L^2$ norm interpolation. The duality principle provides a natural and constructive avenue for lower bounding their singular values. Since we have no exploitable information on the singular vectors of either operator, we construct interpolants for arbitrary data $v$ on $\calX$, and then estimate these interpolants in $L^2$ uniformly in $v$. This leads us to the subsequent definition and lemma. 

\begin{definition}
	For any set $\calX=\{x_k\}_{k=1}^s\subset [-\frac 12,\frac 12)^d$, we say $\{f_k\}_{k=1}^s$ is a family of Lagrange interpolants for $\calX$ if $f_k(x_\ell)=1$ if $k=\ell$ and $f_k(x_\ell)=0$ if $k\not=\ell$. 
\end{definition}

\begin{lemma}
	\label{lem:duality2}
	Suppose $\Omega\subset\R^d$ is compact and $\calX\subset [-\frac 12,\frac 12)^d$ has cardinality $s$. If there is a family $\{f_k\}_{k=1}^s\subset \calW(\Omega)$ of Lagrange interpolants for $\calX$, then 
	$$
	\sigma_{\min}(\calF_{\Omega,\calX})
	\geq \frac 1 {\sqrt s} \, \min_{1\leq k\leq s} \|f_k\|_{L^2(\R^d)}^{-1}.
	$$
	If there is a family $\{f_k\}_{k=1}^s\subset \calP(\Omega)$ of Lagrange interpolants for $\calX$, then
	$$
	\sigma_{\min}(\Phi_{\Omega,\calX})
	\geq \frac 1 {\sqrt s} \, \min_{1\leq k\leq s} \|f_k\|_{L^2(\T^d)}^{-1}.
	$$
\end{lemma}

This lemma is proved in \cref{proof:duality2} and it allows us to reduce the problem of lower bounding their singular values to constructing Lagrange interpolants with suitably controlled norms. We briefly mention a simple observation that will allow us to handle the two operators in a unified manner. Notice that if we extend $f\in \calP(\Omega)$ periodically to $\R^d$, then $f\in L^\infty(\R^d)$ but $f\not\in L^2(\R^d)$. So we can only take its Euclidean Fourier transform in the sense of distributions. Nonetheless, if $g\in \calW(\Omega_2)$, then 
$$
fg\in L^2(\R^d) \andspace fg\in \calW((\Omega_1\cap \Z^d)+\Omega_2),
$$
where $+$ is the algebraic sum of two sets. Due to this observation, we will mainly focus on constructing trigonometric polynomial interpolants and slightly modify them to get bandlimited ones.  

The next proposition acts as a converse to \cref{prop:duality}, and shows that any lower bound on the smallest singular value provides the existence of polynomials with prescribed interpolation properties. A proof of the following lemma for $d=1$ can be found in \cite[Proposition 5.3]{li2023multiscale} and readily extends to higher dimensions. 

\begin{lemma}
	\label{prop:interpolation2}
	Suppose $\Omega\subset\R^d$ is compact, $\calX\subset [-\frac 12,\frac 12)^d$ has cardinality $s$, and $\sigma_{\min}(\Phi_{\Omega,\calX})>0$. Then for any $w\in\C^s$, there exits $f\in \calP(\Omega)$ such that $f|_\calX=w$, 
	$$
	\|f\|_{L^2(\T^d)}\leq \frac{|w|_2}{\sigma_{\min}(\Phi_{\Omega,\calX})} \andspace 
	\|f\|_{L^\infty(\T^d)}\leq \sqrt{|\Omega|_*} \, \frac{|w|_2}{\sigma_{\min}(\Phi_{\Omega,\calX})}.
	$$ 
\end{lemma}

We loosely refer to the strategy provided by the above results as the {\it polynomial (interpolation) method}. A primary usefulness of this connection between Fourier matrices and trigonometric interpolation is that it enables us to use tools from Fourier analysis and polynomial approximation, instead of solely working with matrices. While the polynomial method is conceptually helpful, it is only useful if we can construct Lagrange interpolants with small $L^2$ norm, otherwise the resulting lower bounds for $\sigma_{\min}(\Phi_{\Omega,\calX})$ and $\sigma_{\min}(\calF_{\Omega,\calX})$ would have loose constants and/or suboptimal rates. To carry out explicit constructions, we start with the following decomposition. 

\begin{proposition}
	\label{prop:decomp}
	For any $\tau>0$, $p\in [1,\infty]$, and non-empty $\calX\subset [-\frac 12,\frac 12)^d$, letting $\nu:=\nu_p(\tau,\calW)$, there exist non-empty disjoint subsets $\calX_1,\calX_2,\dots,\calX_\nu\subset\calX$ such that their union is $\calX$ and $\Delta_p(\calX_k)>\tau$ for each $k=1,\dots,\nu$.
\end{proposition}

This proposition can be readily deduced from a greedy construction of the $\calX_1,\dots,\calX_\nu$, see the proof of \cite[Proposition 6.1]{li2023multiscale} for $d=1$. The next proposition illustrates a nontrivial application of the polynomial method and is proved in \cref{proof:localization}. 

\begin{lemma}[Localization]
	\label{prop:localization}
	Let $\calX:=\{x_k\}_{k=1}^s\subset [-\frac 12,\frac 12)^d$, $m\in \N_+$, and $p\in [1,\infty]$. Assume there exist $C>0$, $c\in (0,1)$, and $\tau>0$ such that 
	\begin{equation}
		\label{eq:densitycondition}
		\frac{C \nu_p(\tau,\calX)}{\tau} \leq \frac m 2,
	\end{equation}	
	and for any $\calW\subset \T^d$ with $\Delta_p(\calW)\geq \tau$, we have 
	\begin{equation}
		\label{eq:wellsepbound}
		\sigma_{\min}\left(\Phi_{\Omega^p_{C/\tau}, \, \calW}\right)
		\geq c\sqrt{ |\Omega^p_{C/\tau}|_*}. 
	\end{equation}
	For each $k\in \{1,\dots,s\}$, there exists $g_k\in \calP(\Omega^p_{m/2})$ such that $g_k(x_k)=1$, $g_k$ vanishes on $\calX\setminus \calN_p(x_k,\tau,\calX)$, and
	$$
	\|g_k\|_{L^\infty(\T^d)} \leq \frac 1 {c^{\nu_p(\tau,\calX)}}. 
	$$
	Consequently, if $\sigma_{\min}(\Phi_{\Omega_{m/2}^p,\calN_p(x_k,\tau,\calX)})>0$ for each $k\in \{1,\dots,s\}$, then we have
	\begin{equation}
		\label{eq:localization}
		\sigma_{\min}(\Phi_{\Omega_m^p,\calX})
		\geq \frac {c^{\nu_p(\tau,\calX)}} {\sqrt s} \min_{1\leq k\leq s}  \sigma_{\min}\big(\Phi_{\Omega^p_{m/2},\calN_p(x_k,\tau,\calX)}\big). 
	\end{equation}
\end{lemma}

This lemma shows that if one has a lower bound for the smallest singular value in the well separated case \eqref{eq:wellsepbound}, then together with an assumption on the local sparsity \eqref{eq:densitycondition}, one gets polynomials $\{g_k\}_{k=1}^s$ that vanish on all points in $\calX$ that are at least distance $\tau$ away from $x_k$. These polynomials are useful because they allow us to reduce the global problem of estimating  $\sigma_{\min}(\Phi_{\Omega^p_m,\calX})$ to a collection of $s$ local problems, one for each neighborhood set $\calN_p(x_k,\tau,\calX)$. 

Provided that the assumptions hold, this localization comes at a price of course. Comparing both sides of the inequality in \eqref{eq:localization}, the number of Fourier samples are reduced from $|\Omega^p_m|_*$ to $|\Omega^p_{m/2}|_*$, so resolution worsens by a factor of 2. We also pick up an additional $c^{\nu_p(\tau,\calX)}$ term, which can be small if either $c$ is close to zero and/or $\nu_p(\tau,\calX)$ is large. For many sets of interest, such as clumps, the $\tau$ local sparsity of $\calX$ is small. When we apply this theorem, the constant $c$ will come from \cref{thm:wellsepball,thm:wellsepcube}, so it is important to obtain good constants for the well separated case to temper the effects of localization.

With the localization lemma at hand, let us explain at a high level how to deal with the local set $\calN_p(x_k,\tau,\calX)$. We will again employ the duality principle, which transforms the problem to constructing Lagrange polynomials $\{f_k\}_{k=1}^s$ for $\calN_p(x_k,\tau,\calX)$. Shifting by $x_k$ and noting that both $\calP(\Omega)$ and $\calW(\Omega)$ are invariant under shifts, this reduces the problem down to the canonical case whereby $0\in \calU$ and all other elements in $\calU$ are sufficiently close to zero. Importantly, we place no lower bounds on $\Delta_p(\calU)$. We seek a $f$ such that $f(0)=1$, $f$ vanishes on $\calU\setminus \{0\}$, and its $L^2$ norm is suitably controlled. 

Let us briefly describe our strategy for constructing a reasonable interpolant $f$.  Given a nonzero $u\in \calU$, a natural strategy is to find a suitable $a:=a(u)\in\R^d$ and let
\begin{equation}
	\label{eq:basicinter}
	f(x):=\frac{e^{2\pi i a\cdot x}-e^{2\pi i a\cdot u}}{1-e^{2\pi i a\cdot u}}.
\end{equation}
This interpolates the data points $(0,1),(u,0)\subset [-\frac 12,\frac 12)^d\times \R$. Doing this for each $u\in \calU\setminus\{0\}$ and multiplying them together yields a Lagrange interpolant. 

\begin{figure}[h]
	\centering
	\begin{tikzpicture}[scale=.5]
		\filldraw (0,0) circle (3pt) node[anchor=north east]{$0$};
		\filldraw (4,2) circle (3pt) node[anchor=north east]{$u$};
		\draw[thick,->] (4,2)--(8,4) node[anchor=south east]{$a$};
		\draw[thick,->,dashed] (4,2)--(8,3.6) node[anchor=north]{$q$};
		\draw[thick] (-2,4)--(2,-4) node[anchor=north west]{$f=1$};
		\draw[thick] (2,6)--(6,-2) node[anchor=north west]{$f=0$};
		\draw[thick,dashed] (6,-3)--(2,7) node[anchor=south]{$\tilde f=0$};
		\draw[thick,dashed] (2,-5)--(-2,5) node[anchor=south]{$\tilde f=1$};
	\end{tikzpicture}
	\caption{The sets $\{f=1\}$ and $\{f=0\}$ shown in solid are orthogonal to $a$, while the sets $\{\tilde f=1\}$ and $\{\tilde f=0\}$ shown in dashed are orthogonal to $q$.}
	\label{fig:localinter}
\end{figure}
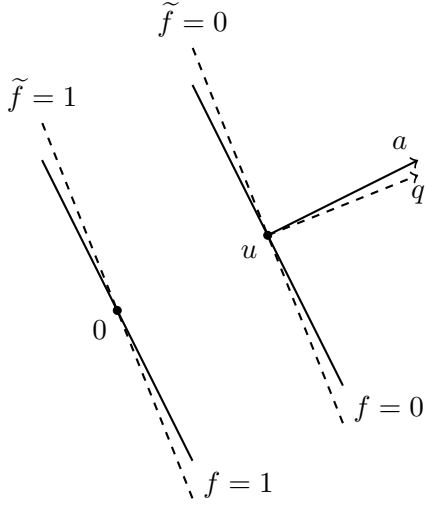

This construction would be quite straightforward and would be analogous to the one-dimensional case, but there are two issues. First, $f$ is not a trigonometric polynomial if $a(u)\not\in\Z^d$, which is problematic in view of the duality principle for $\Phi_{\Omega,\calX}$. This creates an additional complication of replacing the ``analog" $a(u)$ with a suitable ``quantized" $q(u)\in \Z^d$. Hence, the quantized function is 
\begin{equation}
	\label{eq:basicinter2}
	\tilde f(x):=\frac{e^{2\pi i q\cdot x}-e^{2\pi i q\cdot u}}{1-e^{2\pi i q\cdot u}}.
\end{equation}
Second, we need to pick $q(u)$ carefully so that $|q(u)\cdot u|$ is not too close to zero, in order to ensure that $\|\tilde f\|_{L^\infty(\T^d)}$ is reasonably controlled, yet $|q(u)|_p$ itself is not too large otherwise the obtained interpolant might not be in $\calP(\Omega)$. See Figure \ref{fig:localinter} for a cartoon of how replacing $a$ with $q$ effectively changes the direction of $f$.

\begin{lemma}
	\label{lem:quanerror}
	Let $d\geq 2$ be an integer, $p\in [1,\infty]$, and $u\in [-\frac 12, \frac 12)^d$ with $0<|u|_{p'}\leq \frac {1}{4d^{1/p}}$. For any $\alpha>0$ such that $|u|_{p'}\leq \alpha \leq \frac {1}{4d^{1/p}}$, there exists $q\in \Z^d$ such that 
	\begin{equation}
		\label{eq:lemquanerror} 
		|q|_p\leq \frac{1}{2\alpha}, \quad 
		\frac 1 {4\alpha} |u|_{p'} \leq |q\cdot u| \leq \frac 1 2, \andspace 
		|1-e^{2\pi i q \cdot u}| \geq \frac{\sqrt 2}{ \alpha} |u|_{p'}.  
	\end{equation}
\end{lemma}

With this lemma at hand, we will first consider a general situation where no additional information about $\calU$ is provided. Consequently, the upper bound for $\|f\|_{L^\infty(\T^d)}$ that we will obtain naturally corresponds to worst cases of $\calU$. 

\begin{lemma}
	\label{lem:neighborset1}
	Let $p\in [1,\infty]$ and $\calU\subset [-\frac 12, \frac 12)^d$ be a finite set of at most $r$ elements such that $0\in \calU$ and $|u|_{p'}\leq \frac 1{4 d^{1/p}}$ for each $u\in \calU$. For any real $n\geq 2 d^{1/p}r$, there exists $f\in \calP(\Omega^p_{n(r-1)/r})$ such that $f(0)=1$, $f$ vanishes on $\calU\setminus \{0\}$, and
	\begin{equation}
		\label{eq:neighbor1}
		\|f\|_{L^\infty(\T^d)}
		\leq \sqrt{2^{|\calU|-1}} \prod_{0<|u|_{p'} \leq \frac {r} {2n}} \, \frac {r} {2n |u|_{p'}}.
	\end{equation}
\end{lemma}

Its proof can be found in \cref{proof:neighborset1}. Now we are ready to state our main theorems for the super-resolution regime when Fourier samples are collected in a ball or cube. These theorems are proved in \cref{proof:srball,proof:srcube} respectively. Recall the quantity $c(\alpha)$ in \eqref{eq:Calphadef} and \eqref{eq:Calphadef2}. 

\begin{theorem}[Super-resolution regime for arbitrary nodes and samples in a ball]	\label{thm:srball} 
	Let $d\geq 2$, $\alpha>\frac{j_{d/2-1,1}}\pi$, and $\calX:=\{x_k\}_{k=1}^s\subset [-\frac 12,\frac 12)^d$. Suppose $m\geq 4 s\sqrt d$ and there is a $0<\tau\leq \frac 1{4\sqrt d}$ such that $\frac{2 \alpha \nu}{\tau} \leq m$, where $\nu:=\nu_2(\tau,\calX).$	Then we have the inequalities,
	\begin{align*}
		\sigma_{\min}(\calF_{B_m,\calX})
		&\geq \frac 1{\sqrt{2^{\nu-1} s}}\(\frac{c(\alpha)}{|B_1|}\)^{\nu/2} \(\frac{|B_{\alpha/\tau}|}{|B_{\alpha/\tau}|_*}\)^{\nu/2} \, \sqrt{|B_{m/(2\nu)}|} \, \,  \min_{1\leq k\leq s} \Bigg\{\prod_{0<|x_j-x_k|_2\leq \frac{\nu}{m}} \, \frac {m}{\nu} \, |x_j-x_k|_2\Bigg\}, \\
		\sigma_{\min}(\Phi_{B_m,\calX})
		&\geq \frac 1{\sqrt{2^{\nu-1} s}}\(\frac{c(\alpha)}{|B_1|}\)^{\nu/2} \(\frac{|B_{\alpha/\tau}|}{|B_{\alpha/\tau}|_*}\)^{\nu/2} \sqrt{|B_{m/(2\nu)}|_*} \, \min_{1\leq k\leq s} \Bigg\{\prod_{0<|x_j-x_k|_2\leq \frac{\nu}{m}} \, \frac {m}{\nu} \, |x_j-x_k|_2\Bigg\}.
	\end{align*}
\end{theorem}

\begin{theorem}[Super-resolution regime for arbitrary nodes and samples in a cube]
	\label{thm:srcube}
	Let $d\geq 2$ and $\beta\geq \frac 1 {2 \log 2}$, and $\calX:=\{x_k\}_{k=1}^s\subset [-\frac 12,\frac 12)^d$. Suppose $m\geq 4 s$ and there is a $0<\tau\leq \frac 1 {4d}$ such that $\frac{2\beta d\nu}{\tau} \leq m$, where $\nu:=\nu_\infty(\tau,\calX)$. Then we have the inequalities,
	\begin{align*}
		\sigma_{\min}(\calF_{Q_m,\calX})
		&\geq \sqrt{\frac 2 s} \, \( 1-\frac 12 e^{1/(2\beta)}\)^{\nu/2} \, \sqrt{|Q_{m/(2\nu)}|} \, \,  \min_{1\leq k\leq s} \Bigg\{\prod_{0<|x_j-x_k|_1\leq \frac{\nu}{m}} \, \frac {m}{\nu} \, |x_j-x_k|_1  \Bigg\}, \\
		\sigma_{\min}(\Phi_{Q_m,\calX})
		&\geq \sqrt{\frac 2 s} \, \( 1-\frac 12 e^{1/(2\beta)}\)^{\nu/2} \sqrt{|Q_{m/(2\nu)}|_*} \, \min_{1\leq k\leq s} \Bigg\{\prod_{0<|x_j-x_k|_1\leq \frac{\nu}{m}} \, \frac {m}{\nu} \, |x_j-x_k|_1  \Bigg\}.
	\end{align*}
\end{theorem}

We first discuss their similarities. The terms $\sqrt{|B_m|}$ and $\sqrt{|Q_m|}$ are natural scaling terms, which are proportional to $\sqrt{|B_{m/(2\nu)}|}$ and $\sqrt{|Q_{m/(2\nu)}|}$ modulo constants that depend on $\nu$ and $d$ but not on $m$. Notice that each term in theses products is at most 1, so their smallest singular values are both small if there are many elements of $\calX$ that are close together. For example, if for each $k\in \{1,\dots,s\}$ we define 
$s_{k,\ell}$ to be the number of elements in $\calX\setminus \{x_k\}$ whose $\ell^p$ distance to $x_k$ is in the interval $(\frac 1 {2^{\ell+1}} \frac \nu m,\frac 1 {2^\ell}\frac \nu m]$, then
$$
\log_2 \Bigg( \prod_{0<|x_j-x_k|_p\leq \frac{\nu}{m}} \, \frac {m}{\nu} \, |x_j-x_k|_p \Bigg)
\asymp -\sum_{j=0}^\infty s_{k,\ell} \, \ell.
$$ 
This shows that interactions between $\calX$ at small scales dominate the product term. This way of writing it may be more familiar to inequalities that appear in multiscale analysis (e.g., wavelet or Littlewood-Paley characterizations of function spaces). We call these theorems multiscale estimates precisely for this reason. They are higher dimensional generalizations of the multiscale estimates given in \cite{li2023multiscale} for the one-dimensional setting.

Now, we discuss the main differences between \cref{thm:srball,thm:srcube}. The primary difference is the metric. When Fourier samples are collected in a ball, all distances are measured in  $\ell^2$, which is natural since the ball $B_m$ and the $\ell^2$ distance are rotationally invariant. However, when Fourier samples are collected in a cube, the local sparsity $\nu$ is defined using the $\ell^\infty$ distance, while all small scales are measured in $\ell^1$. This is natural since $\ell^1$ and $\ell^\infty$ are dual. It also shows that taking measurements in a cube offers an advantage over the ball because $|x|_2\leq |x|_1$, so elements in $\calX$ are more separated in $\ell^1$ as opposed to $\ell^2$. There are some slight differences in the requirements for $m$ and $\tau$, which are equivalent up to constants that depend on $p$ and $d$. 

For the ball case, the constants that appear before the product terms are more complicated. Results for the generalized Gauss circle problem tell us that, for fixed dimension $d\geq 2$, the leading order growth rates of $|B_r|$ and $|B_r|_*$ are identical in the limit $r\to\infty$. See \cite{ivic2004lattice} for a survey on this topic. Consequently, $|B_{\alpha/\tau}|/|B_{\alpha/\tau}|_*$ can be approximated by 1 when $\alpha/\tau$ is large. There are numerous open and difficult questions about how tight this approximation is. In several examples that we discuss later on, $\alpha/\tau$ can be selected proportional to $m$, so this quotient term can be effectively ignored when $m$ is large. The constant $c(\alpha)$ is small when $d$ is large. This can be addressed using an alternative result (that is not included in this paper) for the well separated case for samples in a ball, see the discussion following \cref{thm:wellsepball}.

An important aspect of both theorems is the choice of $\tau$. In general, there may be infinitely many valid choices, and ideally, one selects $\tau$ optimally to maximize the provided lower bounds. To do this, one would need to search through all possible $\tau$ and evaluate these expressions and take the max. 

For certain sets, we can select $\tau$ heuristically. Clumps are natural examples of sets that satisfy the assumptions of our main theorems. Note that if $\calX$ consists of clumps with parameters $(p,\tau,\lambda)$, then $\nu_p(\tau,\calX)= \lambda$ because any two clumps are at least $\tau$ apart and at least one clump has cardinality $\lambda$. If we also include an extraneous assumption that the minimum separation is bounded by $\frac \lambda m$, then the product terms that appear in \cref{thm:srball,thm:srcube} can be lower bounded by $(\frac{m\delta}\lambda)^{\lambda-1}$. From these observations, we immediately get the following corollaries.

\begin{corollary}[Super-resolution regime for clumps and samples in a ball]
	\label{cor:clumpsball}
	Let $d\geq 2$ and $\alpha>\frac{j_{d/2-1,1}}\pi$. Suppose $\calX\subset\T^d$ has cardinality $s$ and consists of clumps with parameters $(2,\tau,\lambda)$. Suppose $m\geq 4s\sqrt d$, $\Delta_2(\calX)=:\delta\leq \frac \lambda m$, and $\frac{2\alpha \lambda}{m}\leq \tau\leq \frac 1 {4 \sqrt d}$. Then we have 
	\begin{align*}
		\sigma_{\min}(\calF_{B_m,\calX})
		&\geq \frac 1{\sqrt{2^{\lambda-1}s}}\(\frac{c(\alpha)}{|B_1|}\)^{\lambda/2} \(\frac{|B_{\alpha/\tau}|}{|B_{\alpha/\tau}|_*}\)^{\lambda/2} \, \sqrt{|B_{m/(2\lambda)}|} \, \, \,  \( \frac{m\delta}{\lambda}\)^{\lambda-1}, \\
		\sigma_{\min}(\Phi_{B_m,\calX})
		&\geq \frac 1{\sqrt{2^{\lambda-1}s}}\(\frac{c(\alpha)}{|B_1|}\)^{\lambda/2} \(\frac{|B_{\alpha/\tau}|}{|B_{\alpha/\tau}|_*}\)^{\lambda/2} \sqrt{|B_{m/(2\lambda)}|_*} \, \, \( \frac{m\delta}{\lambda}\)^{\lambda-1} .
	\end{align*}	
\end{corollary}

\begin{corollary}[Super-resolution regime for clumps and samples in a cube]
	\label{cor:clumpscube}
	Let $d\geq 2$ and $\beta\geq \frac 1 {2 \log 2}$. Suppose $\calX\subset\T^d$ has cardinality $s$ and consists of clumps with parameters $(\infty,\tau,\lambda)$. Suppose $m\geq 4s$, $\Delta_1(\calX)=:\delta\leq \frac \lambda m$, and $\frac{2\beta d\lambda}{m}\leq \tau\leq \frac 1 {4d}$. Then we have  
	\begin{align*}
		\sigma_{\min}(\calF_{Q_m,\calX})
		&\geq \sqrt{\frac 2 s} \, \( 1-\frac 12 e^{1/(2\beta)}\)^{\lambda/2} \, \sqrt{|Q_{m/(2\lambda)}|} \, \, \,  \( \frac{m\delta}{\lambda}\)^{\lambda-1}, \\
		\sigma_{\min}(\Phi_{Q_m,\calX})
		&\geq \sqrt{\frac 2 s} \, \( 1-\frac 12 e^{1/(2\beta)}\)^{\lambda/2} \, \sqrt{|Q_{m/(2\lambda)}|_*} \, \, \( \frac{m\delta}{\lambda}\)^{\lambda-1} .
	\end{align*}	
\end{corollary}

The above dependence on $m$ and $\delta$ are optimal. For example, suppose the elements of $\calX$ are equally spaced by $\delta$ on a line parallel to one of the coordinate axes. Then the matrix $\Phi_{Q_m,\calX}$ consists of $(2m+1)^{d-1}$ copies of a one-dimensional Fourier matrix corresponding to nodes $\{0,\delta,\dots,(\lambda-1)\delta\}$ and frequencies $\{-m,\dots,m\}$. Upper bounds for the smallest singular value of one-dimensional Fourier matrices, see \cite[Proposition 2.10]{li2021stable} tell us that the smallest singular value of each one is upper bounded by $C_\lambda \sqrt m (m\delta)^{\lambda-1}$. Using this for each of the $(2m+1)^{d-1}$ repeated matrices yields the upper bound 
$$
\sigma_{\min}(\Phi_{Q_m,\calX})
\lesssim_{\lambda,d} \sqrt{|Q_m|} (m\delta)^{\lambda-1}. 
$$
Similar observations hold for the ball. Note that $B_m\subset Q_m$ implies $\sigma_{\min}(\Phi_{B_m,\calX})\leq \sigma_{\min}(\Phi_{Q_m,\calX})$. Also $|B_m|$ and $|Q_m|$ are equivalent up to a constant that only depends on $d$, so we have
$$
\sigma_{\min}(\Phi_{B_m,\calX})
\lesssim_{\lambda,d} \sqrt{|B_m|} (m\delta)^{\lambda-1}. 
$$

A similar result to \cref{cor:clumpscube} was established in Kunis-Nagal \cite{kunis2020smallest}. They lower bounded $\sigma_{\min}(\Phi_{Q_m',\calX})$ where $Q_m'=\{0,\dots,m\}^d$ and obtained the same $C_{\lambda,d} m^{d/2}(m\delta)^{\lambda-1}$ lower bound up (with different implicit constant), but under significantly stronger assumptions on the clump separation $\tau$. Indeed, their weakest condition on $\tau$, see \cite[Examples 4.3 (ii)]{kunis2020smallest}, requires 
$$
\tau \geq \frac {3.3 \lambda ( 2.5+1.4d + \lambda \log \lambda+ \log \mathfrak{C})}{m}, \wherespace
\mathfrak{C} := \max_{1\leq k\leq s} \Bigg\{\prod_{0<|x_j-x_k|_\infty\leq \frac{1}{m}} \, \frac {1}{m|x_j-x_k|_\infty}  \Bigg\}.
$$ 
As the minimum separation goes to zero, necessarily $\mathfrak{C}\to \infty$, which implies $\tau\to\infty$. So this condition becomes vacuous if any two nodes in $\calX$ are pushed closer together, while the rest remain fixed. This is a proof artifact that can be traced back to the interpolation techniques developed in \cite{li2021stable}. To compare, 
the separation condition in \cref{cor:clumpscube} does not depend on $\delta$, scales linearly with $\lambda$ instead of quadratically, and has the same dependence on $d$. Another main difference is that they use the $\ell^\infty$ distance on $\calX$, whereas our corollary uses the $\ell^1$ metric, which is preferable since $|x|_\infty\leq |x|_1$. Our improved results for multidimensional clumps are largely due to new technical developments for the one-dimensional problem provided by \cite{li2023multiscale}, which dramatically improved upon \cite{li2021stable}. 

During the final stages of writing this draft, we became aware that Diab-Batenkov \cite{diab2024spectral} are working on a related topic. They focus on an arbitrary fixed $\calX$ and examine all the singular values of $\Phi_{Q_m,\calX_\delta}$ as $\delta\to0$, where $\calX_\delta:=\{\delta x\colon x\in\calX\}$ is an isotropic dilation of $\calX$. Since $\delta\to0$, this is a single clump configuration corresponding to $\lambda=s$. They compute the number of singular values that decay at a rate $\delta^k$. For the smallest singular value, they obtain the rate $\delta^{s-1}$, which matches our result when $\calX$ is a single clump. The main difference is that we obtain non-asymptotic bounds that pertain to more general configurations, but we only have results for the smallest singular value. 

\section{Super-resolution regime and hyperplane geometry}
\label{sec:sr2}

The material contained in this section is motivated  by the following experiment and numerical observations. Fix any $m\in \N_+$. As we know from previous results, the local behavior of $\calX$ dominates the smallest singular values of both Fourier operators. For simplicity, suppose $\calU\subset Q_{1/m}$ and has cardinality $\lambda$, so in the context of this paper, $\calU$ represents a local subset of the true support. We define the isotropic dilation of $\calU$ by $\delta\in (0,1)$ as  
\begin{equation}
	\label{eq:Xgeneric}
	\calU_\delta:=\{\delta x\colon x \in\calU\}. 
\end{equation}
In the previous section, we know that for worst case $\calU$, we have $\sigma_{\min}(\Phi_{Q_m,\calU_\delta})\asymp \delta^{\lambda-1}$ as $\delta\to 0$. The main question is whether this rate is sharp for typical or generic $\calU$. 

\begin{table}[h]
	\centering
	\begin{tabular}{|l|c|}\hline
		$\lambda$ &$\gamma(\lambda,2)$ \\ \hline 
		1 &0\\ \hline 	
		2, 3 &1 \\ \hline 	
		4, 5, 6 &2 \\ \hline 	
		7, 8, 9, 10 &3 \\ \hline 	
		11, 12, 13, 14, 15 &4 \\ \hline 	
		16, 17, 18, 19, 20, 21 &5 \\ \hline 		
	\end{tabular} \quad 
	\begin{tabular}{|l|c|}\hline
		$\lambda$ &$\gamma(\lambda,3)$ \\ \hline 
		1 &0 \\ \hline 	
		2, 3, 4 &1 \\ \hline 	
		5, 6, 7, 8, 9, 10 &2 \\ \hline 	
		11, 12, \dots, 19, 20 &3 \\ \hline 	
		21, 22, \dots, 34, 35 &4 \\ \hline 	
		36, 37, \dots, 55, 56  &5\\ \hline 		
	\end{tabular}
	\caption{Computed powers $\gamma(\lambda,d)$. Left is $d=2$ and right is $d=3$.}
	\label{table:exponent}
\end{table}

To carry out a motivational experiment, we let $\calU$ consist of $\lambda$ independent draws from the uniform distribution on $Q_{1/m}$ for dimensions $d\in \{2,3\}$. As we vary $\delta$ for each $(\lambda,d)$, we numerically see that for an exponent $\gamma(\lambda,d)$ listed in Table \ref{table:exponent}, that 
\begin{equation}
	\label{eq:randrate}
	\sigma_{\min}(\Phi_{Q_m,\calU_\delta})\asymp \delta^{\gamma(\lambda,d)}.
\end{equation}
It is important to mention that the reported $\gamma(\lambda,d)$ is consistent across random draws of $\calU$ -- for each trial, we obtained the exact same exponents and only the implicit constant in front of $\delta^{\gamma(\lambda,d)}$ depended on the realization $\calU$. We also obtained the same exponents when Fourier samples are drawn from $B_m$ instead of $Q_m$, and when $\calU$ are independent samples from other absolutely continuous distributions such as a normal distribution. These results indicate that \eqref{eq:randrate} is not a concentration phenomenon and is a manifestation of some property that all typical sets enjoy. The numerical results here emphatically show that isotropic dilations of a typical set does not achieve the worst case rate of $\delta^{\lambda-1}$, especially for large $\lambda$. 

After completing this experiment, we found out that the asymptotic rate of $\delta^{\gamma(\lambda,d)}$ for dilations of a generic set was rigorously established by Barthelm\'e-Usevich \cite[Theorem 6.1]{barthelme2021spectral}. They showed that $\gamma(\lambda,d)$ is precisely the smallest integer $\gamma$ such that ${\gamma+d \choose d }\geq \lambda$. This is consistent with the numerically computed exponents in \cref{table:exponent}. See \cite{diab2024spectral} for some extensions to dilations of arbitrary sets.

We give a different perspective on this phenomenon based on hyperplane decompositions. Our results will be applicable to more general sets beyond $\calU_\delta$ that was used in this experiment, and we will provide explicit non-asymptotic lower bounds. Specializing to dilations of a generic set, our results will provide a different lower bound,
$$
\sigma_{\min}(\calF_{Q_m,\calU_\delta})\gtrsim \delta^{r(\lambda,d)}, \wherespace r(\lambda,d):=\Big\lceil \frac{\lambda-1} d \Big\rceil.
$$
The exponents $\gamma(\lambda,d)$ and $r(\lambda,d)$ are equal for small enough $\lambda$ depending on $d$, but become drastically different as $\lambda$ increases. For instance, $\gamma(\lambda,2)=r(\lambda,2)$ for $\lambda \in\{1,\dots,5,7\}$ and $\gamma(\lambda,3)=r(\lambda,3)$ for $\lambda \in \{1,\dots,7\}$. For imaging applications, $\lambda$ represents the number of points within a $\ell^p$ ball whose radius is on the order of a Rayleigh length $\frac 1{2m}$, so small $\lambda$ such as $\lambda=2,3,4$ is already challenging for super-resolution imaging. 

The quantity $r(\lambda,d)$ is related to hyperplanes. Again, looking at this problem from the dual perspective provided by \cref{prop:duality}, a plane wave $x\mapsto e^{2\pi i\omega \cdot x}$ is a canonical ``atom" for $\calW(\Omega)$. It is constant on hyperplanes and any set of $d$ generic points in $\R^d$ specify a unique hyperplane that contains them. Hence, if $\calX\setminus\{x_k\}$ consist of at most $\lambda-1$ elements near $x_k$, we can find $r(\lambda,d)$ hyperplanes that contain them. Importantly, if $\calX$ is locally a generic set, then these hyperplanes do not intersect $x_k$. So to build a local Lagrange interpolant in $\calW(\Omega)$ for a neighborhood of $x_k$, we actually only need to multiply together $r(\lambda,d)$ many plane waves. The following constructs said Lagrange interpolants for sets at small scales and is proved in \cref{proof:neighborset2}.  

\begin{lemma}
	\label{lem:neighborset2}
	Let $d\geq 2$, $p\in [2,\infty]$, and $\calU\subset [-\frac 12, \frac 12)^d$ be a finite set such that $0\in \calU$. Suppose $\calU$ consists of $r$ hyperplanes relative to 0, and let $\eta_1,\dots,\eta_r$ denote their $\ell^2$ distances to 0. For any real $n$ such that $\max_k\eta_k\leq \frac {r+1}{4n}$, there exists $f\in \calW(\Omega_n^p)$ such that $f(0)=1$, $f$ vanishes on $\calU\setminus \{0\}$, and
	\begin{equation}
		\label{eq:neigh2}
		\|f\|_{L^2(\R^d)}
		\leq \sqrt{\frac {2^r}{|\Omega^p_{n/(r+1)}|}} \, \prod_{k=1}^r \frac{r+1}{4 n\eta_k}.
	\end{equation}
\end{lemma}

We make a few comments. First, unlike many of our other results that hold for any $p\in [1,\infty]$, this lemma requires $p\in [2,\infty]$, unless we introduce dimension dependence on $\eta_k$. Second, $\eta_k$ is by definition, the distance between $\calH_k$ and $\{0\}$ in the $\ell^2$ metric. When doing hyperplane decompositions, it seems most natural to work with $\ell^2$ distances regardless of what $\ell^p$ ball the Fourier samples are collected from. It is plausible that this lemma can be adapted for $\ell^q$ hyperplane distances, but we have not attempted to do so. Third, we do not have a corresponding lemma for trigonometric polynomials, which turns out to be much more difficult, and we will return to this issue later. With this lemma at hand, we are in position to state our main results for $\sigma_{\min}(\calF_{\Omega,\calX})$ under hyperplane decompositions of $\calX$. These theorems are proved in \cref{proof:hyper}. Recall the definition of $c(\alpha)$ given in \eqref{eq:Calphadef} and \eqref{eq:Calphadef2}.

\begin{theorem}[Super-resolution regime for local hyperplanes with samples in a ball] \label{thm:hyperball}
	Let $d\geq 2$, $\alpha>\frac{j_{d/2-1,1}}\pi$, $\calX:=\{x_k\}_{k=1}^s\subset [-\frac 12,\frac 12)^d$. Suppose $m\geq 4 s\sqrt d$, there is a $0<\tau\leq \frac 1{4\sqrt d}$ such that 
	$\frac{2 \alpha \nu}{\tau} \leq m$, where $\nu:=\nu_2(\tau,\calX)$, and $\calX$ locally consists of hyperplanes with parameters $(2,\frac \nu m,r,\eta)$ such that $\eta\leq \frac{r+1}{2m}$. Then we have the inequality,
	\begin{align*}
		\sigma_{\min}(\calF_{B_m,\calX})
		&\geq \frac 1{\sqrt{2^{\nu-1}s}}\(\frac{c(\alpha)}{|B_1|}\)^{\nu/2} \(\frac{|B_{\alpha/\tau}|}{|B_{\alpha/\tau}|_*}\)^{\nu/2} \, \sqrt{|B_{m/(2r+2)}|} \, \, \frac{1}{2^{r/2}} \(\frac{2m \eta}{r+1}\)^r.
	\end{align*}
\end{theorem}

\begin{theorem}[Super-resolution regime for local hyperplanes with samples in a cube] \label{thm:hypercube}
	Let $d\geq 2$ and $\beta\geq \frac 1 {2 \log 2}$, and $\calX:=\{x_k\}_{k=1}^s\subset [-\frac 12, \frac 12)^d$. Suppose $m\geq 4 s$, there is a $0<\tau\leq \frac 1 {4d}$ such that $\frac{2\beta d\nu}{\tau} \leq m$, where $\nu:=\nu_\infty(\tau,\calX)$, and $\calX$ locally consists of hyperplanes with parameters $(\infty,\frac \nu m,r,\eta)$ such that $\eta\leq \frac{r+1}{2m}$. Then we have the inequality,
	\begin{align*}
		\sigma_{\min}(\calF_{Q_m,\calX})
		&\geq \sqrt{\frac 2 s} \, \( 1-\frac 12 e^{1/(2\beta)}\)^{\nu/2} \, \sqrt{|Q_{m/(2r+2)}|} \, \, \frac{1}{2^{r/2}} \(\frac{2m \eta}{r+1}\)^r.
	\end{align*}
\end{theorem}

\begin{figure}[ht]
	\centering
	\begin{tikzpicture}[scale=0.4]
		\filldraw (0,0) circle (3pt) node[anchor=north east]{$x_1$};
		\filldraw (4,0) circle (3pt) node[anchor=north east]{$x_2$};
		\filldraw (8,0) circle (3pt) node[anchor=north east]{$x_3$};
		\filldraw (12,0) circle (3pt) node[anchor=north east]{$x_4$};
		\draw[thick] (-1,0)--(13,0);
		\draw[thick,dashed] (4,-2)--(4,2) node[anchor=north west]{$\calH_1$};
		\draw[thick,dashed] (8,-2)--(8,2) node[anchor=north west]{$\calH_2$};
		\draw[thick,dashed] (12,-2)--(12,2) node[anchor=north west]{$\calH_3$};
		\filldraw (25,-1) circle (3pt) node[anchor=north]{$x_1$};
		\filldraw (30,0) circle (3pt) node[anchor=north west]{$x_2$};
		\filldraw (20,0) circle (3pt) node[anchor=north east]{$x_3$};
		\draw[thick] plot[smooth,domain=18:32] (\x, {(\x-25)*(\x-25)/25-1});
		\draw[thick,dashed] (16,0)--(34,0) node[anchor=west]{$\calH_1$};
	\end{tikzpicture}	
	\caption{Points closely spaced on a line and parabola.}
	\label{fig:linepara}
\end{figure}
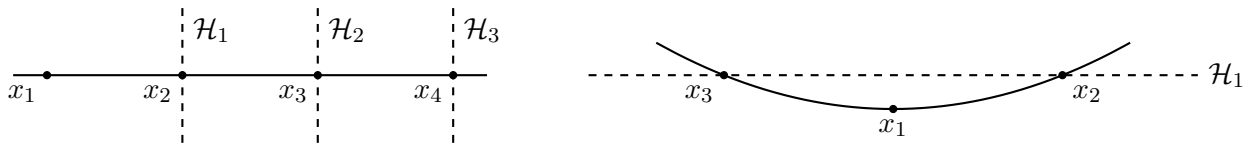

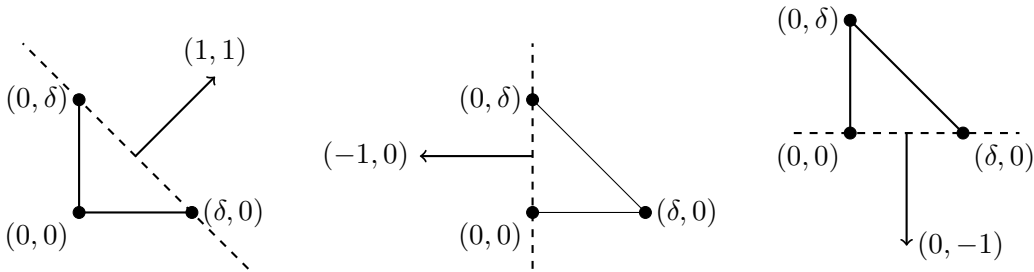
\begin{figure}[b]
	\centering
	\begin{tikzpicture}[scale=.75]
		\filldraw (0,0) circle (3pt) node[anchor=north east]{$(0,0)$};
		\filldraw (2,0) circle (3pt) node[anchor=west]{$(\delta,0)$};
		\filldraw (0,2) circle (3pt) node[anchor=east]{$(0,\delta)$};
		\draw[thick] (0,2)--(0,0)--(2,0); 
		\draw[thick,dashed] (3,-1)--(-1,3);
		\draw[thick,->] (1,1)--(2.41,2.41) node[anchor=south]{$(1,1)$};
	\end{tikzpicture} \quad 
	\begin{tikzpicture}[scale=.75]
		\filldraw (0,0) circle (3pt) node[anchor=north east]{$(0,0)$};
		\filldraw (2,0) circle (3pt) node[anchor=west]{$(\delta,0)$};
		\filldraw (0,2) circle (3pt) node[anchor=east]{$(0,\delta)$};
		\draw (0,2)--(2,0)--(0,0);
		\draw[thick,dashed] (0,-1)--(0,3);
		\draw[thick,->] (0,1)--(-2,1) node[anchor=east]{$(-1,0)$};
	\end{tikzpicture} \quad
	\begin{tikzpicture}[scale=.75]
		\filldraw (0,0) circle (3pt) node[anchor=north east]{$(0,0)$};
		\filldraw (2,0) circle (3pt) node[anchor=north west]{$(\delta,0)$};
		\filldraw (0,2) circle (3pt) node[anchor=east]{$(0,\delta)$};
		\draw[thick,dashed] (-1,0)--(3,0); 
		\draw[thick] (0,0)--(0,2)--(2,0);
		\draw[thick,->] (1,0)--(1,-2) node[anchor=west]{$(0,-1)$};
	\end{tikzpicture}
	\caption{Left to right: the local hyperplane decomposition is chosen with $(0,0)$, $(\delta,0)$ and $(0,\delta)$ as the reference point. A vector with integer coordinates that is orthogonal to each hyperplane is shown as well.}
	\label{fig:triangleex}
\end{figure}

Here we make some basic comparisons between \cref{thm:hyperball,thm:srball} for the ball, and \cref{thm:hypercube,thm:srcube} for the cube. These theorems share similar assumptions and constants that depend on $\tau$, $\nu$, and $\alpha$ or $\beta$, which come from reducing the global problem into to a local one via \cref{prop:localization}. The main difference between these two pairs of theorems is how the local problem is handled. 

\begin{figure}[t]
	\centering
	\begin{subfigure}{0.45\textwidth}					
		\includegraphics[width=\textwidth]{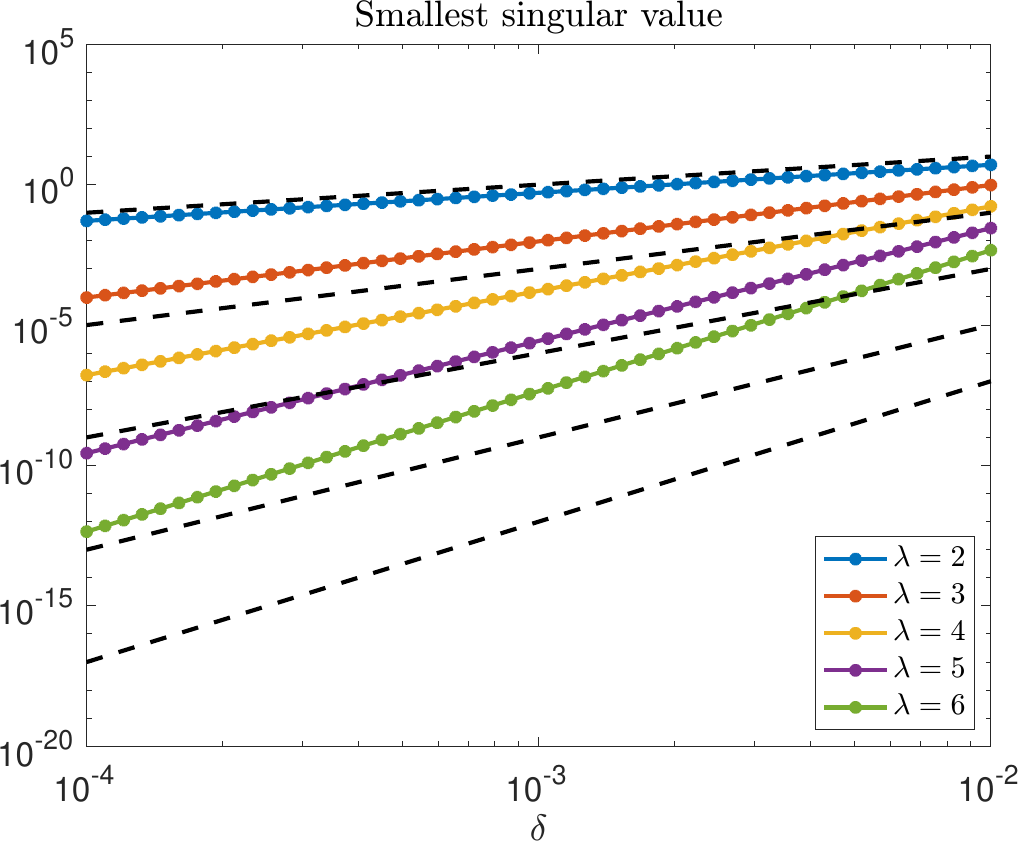}
		\caption{Points on a line}
	\end{subfigure}
	\begin{subfigure}{0.45\textwidth}					
		\includegraphics[width=\textwidth]{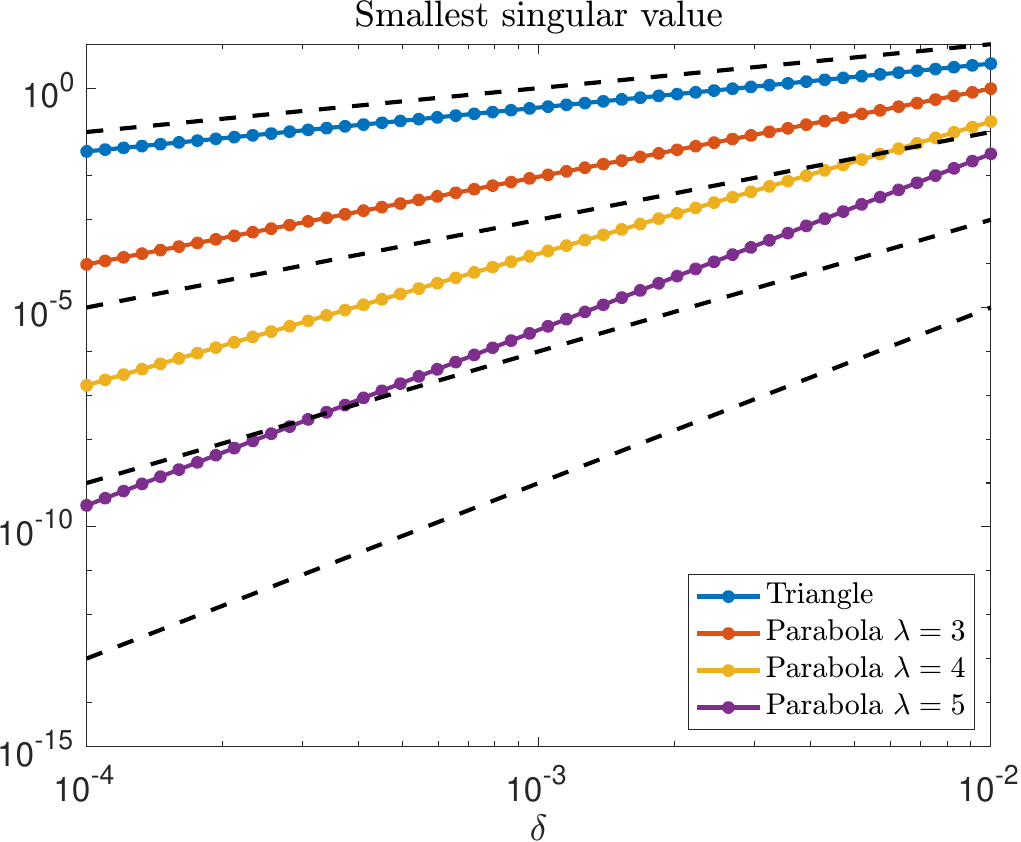}
		\caption{Triangle and parabola}
	\end{subfigure}
	\caption{For $m=20$, plots of $\sigma_{\min}(\Phi_{Q_m,\calX})$ versus scaling factor $\delta$ for several examples of $\calX$, and $C=10^3$. (a) $\calX$ consists of $\lambda$ equally spaced points on a line for $\lambda\in \{2,\dots,6\}$. The dashed lines are $C\delta^k$ for $k\in \{1,\dots,5\}$. (b) $\calX$ is the triangle and parabola examples for $\lambda=2,\dots,5$. The dashed lines are $C\delta^k$ for $k\in \{1,\dots,4\}$. }
	\label{fig:examples}
\end{figure}

Let us look at a few basic examples for $p=\infty$ and $d=2$. They can be adapted to $p=2$ with minor modifications. These examples illustrate some nuances and differences between the theorems in this section versus the previous.
\begin{enumerate}[1.]
	\item 
	Points on a line. Say $\calX$ consists of $\lambda$ points on a line equally spaced by $\delta$. For all sufficiently small $\delta$, we can apply the theorems in \cref{sec:sr1} which predict a rate of $\delta^{\lambda-1}$. The theorems in \cref{sec:sr2} are applicable as well because $\calX$ locally consists of lines. However, we must use at least $r=\lambda-1$ hyperplanes (i.e., lines) since any line that contains two points in $\calX$ will also contain the rest. Hence, we must pick these hyperplanes so that they do not intersect any other points, as shown in Figure \ref{fig:linepara}. Hence, the theorems of this section provide the same $\delta^{\lambda-1}$ rate. The $\delta^{\lambda-1}$ rate is numerically optimal as shown in Figure \ref{fig:examples} (a). 
	
	\item 
	Vertices of a triangle. Say $\calX =\{(0,0),(\delta,0),(0,\delta)\}$ for sufficiently small $\delta$. Here, we see that the theorems in \cref{sec:sr1} predict a rate of $\delta^2$. On the other hand, $\calX$ locally consists of hyperplanes with parameters $(C\delta,p,1,\delta/2)$. Indeed, for any $x\in \calX$, there is a single line that contains $\calX\setminus \{x\}$ and its distance to $x$ is at least $\delta/\sqrt 2$, see Figure \ref{fig:triangleex}. Thus, the theorems in this section give the correct $\delta$ rate. This type of reasoning extends to any arbitrary triangle whose side lengths are on the order of $\delta$. This example was also studied in \cite[Example 5.3]{kunis2020smallest} and the $\delta$ rate for the smallest singular value was derived via explicit calculation, whereas we obtained this using hyperplane methods that can deal with more general sets. 
	
	\item 
	Points on a parabola. Say $\calX =\{(0,0),(-\delta,\delta^2),(\delta,\delta^2)\}$ for sufficiently small $\delta$. Again, the theorems in \cref{sec:sr1} predict a $\delta^2$ rate since the distances between these three points are each on the order of $\delta$. We will see that the hyperplane theorems of this section provide the same results. Indeed, $\calX$ locally consists of hyperplanes with parameters $(C\delta,p,1,\delta^2)$, but $\delta^2$ is the smallest parameter that can be used. See Figure \ref{fig:linepara} where the distance between such a hyperplane and the parabola's vertex scales proportionally to $\delta^2$. Hence, the theorems in this section provide the same rate of $\delta^2$. This example also illustrates that \cref{thm:hyperball,thm:hypercube} and \cref{thm:srball,thm:srcube} may yield the same conclusions even though their exponents are different.
	
	This example can be generalized to $\lambda$ points placed on a sufficiently smooth curve whose distances are proportional to small enough $\delta$. The key observation is that the hyperplane distances can only be lower bounded by $C\delta^2$. For our numerical experiments, we also consider $\lambda=4$ and $\lambda=5$ points on the same parabola $y=x^2$, where the first coordinate of these sets are $-\delta, \, 0,\, \delta, \, 2\delta$ and $-2\delta,\, -\delta,\, 0,\, \delta,\, 2\delta$ respectively. The theorems in this section and the previous one each predict a $\delta^{\lambda-1}$ rate, which is numerically optimal as shown in Figure \ref{fig:examples} (b). This is not surprising since at sufficiently small scales depending on the parabola's curvature, points on a parabola appear as though they are on a line. 
	\item 
	Isotropic dilation of a generic set. Here, we let $\calU_\delta$ be the set in \eqref{eq:Xgeneric}. With probability 1, the set $\calU$ locally consists of hyperplanes with parameters $(C,2,r,c)$, where $r=r(\lambda,2)$ and $C,c>0$ only depend on $\calU$. This implies that the dilation $\calU_\delta$  locally consists of hyperplanes with parameters $(C\delta,2,r,c\delta)$. The theorems of this section predict the rate $\delta^r$, which is significantly better than the $\delta^{\lambda-1}$ rate offered by the theorems in \cref{sec:sr1}. The list of exponents in \cref{table:exponent} generated from numerical simulations indicate that the exponent $r$ is sharp for $\lambda\in \{1,\dots,5,7\}$. 
	
	It is important to mention that neither $\delta^{\gamma(\lambda,d)}$ nor $\delta^{r(\lambda,d)}$ rates hold for a family of anisotropic dilations. For example, three generic points in $\R^2$ will be vertices of a triangle, and if were to dilate them using a diagonal matrix with entries $\delta$ and $\delta^2$, then we would essentially recover the parabola example earlier, which would yield a rate of $\delta^2$.	
\end{enumerate}

Let us discuss some related work. Of course, it is well known that in higher dimensions, the geometry of $\calX$ greatly influences the behavior of polynomial interpolation. Related hyperplane decompositions have been exploited in the context of Lagrange interpolation and appears in books, see \cite[Chapter 10]{cheney2009course} and \cite[Chapter 5]{phillips2003interpolation}. Many prior works focus on finding polynomials of the smallest degree for which interpolation of arbitrary values on $\calX$ is possible. Our focus is in the opposite direction. We allow for large $m$ relative to the local sparsity of $\calX$, which corresponds to using a larger class of functions than necessary for interpolation. In this setting, there is a localization effect of the Fourier transform, so only the local (as opposed to global), hyperplane geometry of $\calX$ predominantly affects $\sigma_{\min}(\calF_{\Omega,\calX})$. 

Somewhat related to the results of this section are analysis of algorithms that aim to recover an unknown measure from its Fourier samples. In the one dimensional setting, the classical Prony's algorithm is a polynomial root finding method \cite{de1795essai}, while methods like MUSIC \cite{schmidt1986multiple}, ESPRIT \cite{roy1989esprit}, and MPM \cite{hua1990matrix} estimate the measure by first identifying an approximation of the range of $\Phi_{\Omega,\calX}$. For both families of algorithms, there are strong performance guarantees for their stability to noise, which can be quantified in terms of one-dimensional Fourier matrices (e.g., Vandermonde with nodes on the unit complex circle), see \cite{moitra2015matrixpencil,li2021stable,li2020super,katz2024accuracy}. While there are multidimensional generalizations \cite{haardt2008higher,kunis2016multivariate,liao2015music} (this list does not include projection based methods that reduce the multidimensional problem into a collection of 1D versions), there are many open questions about their stability to noise, sampling complexity, resolution limits, etc. 

Yet another family of algorithms are those based on convex optimization. They were first introduced in \cite{de2012exact,candes2013super} and later generalized to other measurement operators beyond Fourier and to higher dimensions. Typically, these methods are analyzed through convex duality. One suitable dual formulation is a minimum $L^2$ norm Hermite interpolation problem \cite{duval2015exact}. In contrast, the smallest singular value of both Fourier operators is characterized by a minimum $L^2$ norm Lagrange interpolation problem. These problems are different enough that results for either one does not necessarily carry over to the other.

Despite the differences between these two interpolation problems, there are some common observations. Poon-Peyr\`e \cite{poon2019multidimensional} observed that the configuration of $\calX\subset\R^d$ (i.e., points on a line versus general position) greatly affects the stability of convex algorithms when $\calX$ is dilated by $\delta$ and $\delta\to0$. It is important to remark that they considered general measurement operators beyond Fourier, which could be why they did not observe any relationships to hyperplane decompositions. On the other hand, for Fourier measurements specifically, \cite{benedetto2020super} provided conditions for when solutions to convex algorithms are measures supported in intersections of hyperplanes. 

The rest of this section addresses the natural question of whether these theorems for the continuous Fourier transform also hold for matrices. The numerical simulations in Figure \ref{fig:examples} certainly indicate that this is the case, since they are plots of $\sigma_{\min}(\Phi_{Q_m,\calX})$. To see why this is plausible, both $\calF_{\Omega,\calX}^*\calF_{\Omega,\calX}$ and $\Phi_{\Omega,\calX}^*\Phi_{\Omega,\calX}$ are $s\times s$ matrices and a direct calculation shows that 
\begin{equation}
	\label{eq:FFstar}
	(\calF_{\Omega,\calX}^*\calF_{\Omega,\calX})_{j,k}=\int_{\Omega} e^{2\pi i \omega\cdot (x_j-x_k)} \, d\omega \andspace
	(\Phi_{\Omega,\calX}^*\Phi_{\Omega,\calX})_{j,k}=\sum_{\omega\in \Omega\cap \Z^d} e^{2\pi i \omega\cdot (x_j-x_k)}. 
\end{equation}
For $\Omega=\Omega^p_m$ and large $m$, based off of these formulas, we expect there to be some relationship between their singular values. This perspective yields a simple result for over-sampled Fourier matrices. Note that $\calX\subset [-\frac 12,\frac 12)^d$, so sampling on $\Omega\cap\Z^d$ is natural, but it is a coarse discretization. Fixing an over-sampling factor $\rho\in\N_+$, we define an over-sampled Fourier matrix, 
$$
\Phi_{\Omega,\calX,\rho}:= \Big[ \, \frac 1 {\rho^{d/2}} e^{-2\pi i\omega\cdot x} \, \Big]_{\omega \in \Omega\cap (\rho^{-1}\Z)^d,\, x\in \calX}. 
$$ 
Here, $\rho^{-d/2}$ is the natural normalization for this matrix since there are $\rho^d$ many samples per unit volume cube and the smallest singular value is defined via the $\ell^2$ operator norm. Then a calculation shows that 
\begin{equation}
	\label{eq:PPstar2}
	(\Phi_{\Omega,\calX,\rho}^*\Phi_{\Omega,\calX,\rho})_{j,k} 
	=\frac 1 {\rho^d}\sum_{\omega\in \Omega\cap (\rho^{-1}\Z)^d} e^{2\pi i \omega \cdot (x_j-x_k)}. 
\end{equation}
From here, we see that \eqref{eq:PPstar2} is a Riemann sum approximation of the integral in \eqref{eq:FFstar}, so together with standard matrix perturbation theory, e.g. Weyl's inequality, the singular values of $\Phi_{\Omega,\calX,\rho}$ converge to those of $\calF_{\Omega,\calX}$ as $\rho\to\infty$. This implies that \cref{thm:hyperball,thm:hypercube} readily extend to over-sampled Fourier matrices for large enough $\rho$. 

Obtaining non-asymptotic results for $\sigma_{\min}(\Phi_{\Omega,\calX})$ in terms of hyperplane decompositions poses many technical issues. A primary issue is already evident in the proof of \cref{lem:neighborset2}. In that proof, for each hyperplane, we found a plane wave that vanishes on it, via a Lagrange interpolation formula. It is possible that it does not have integer frequencies (e.g., if a vector orthogonal to the hyperplane has $d-1$ rational entries and 1 irrational entry), so the argument does not hold for trigonometric polynomials. To circumvent this issue, we can place some restrictions on $\calX$. One type of assumption is illustrated by the following lemma, which is proved in \cref{proof:neighborset3}. 

\begin{lemma}
	\label{lem:neighborset3}
	For any integer $d\geq 2$ and $p\in [2,\infty]$, let $\calU\subset [-\frac 12, \frac 12)^d$ be a finite set such that $0\in \calU$. Suppose $\calU\setminus \{0\}$ consists of $r$ hyperplanes $\calH_1\dots,\calH_r$ with parameters $\eta_1,\dots,\eta_r$. Fix any $n$ such that $\max_k\eta_k\leq \frac {r+1}{4n}$ and assume that for each hyperplane $\calH_k$, there is a $q_k\in\Z^d$ orthogonal to $\calH_k$ such that $|q_k|_2\leq \frac {n}{r+1}$. Then there exists $f\in \calP(\Omega_n^p)$ such that $f(0)=1$, $f$ vanishes on $\calU\setminus \{0\}$, and
	\begin{equation}
		\label{eq:neigh3}
		\|f\|_{L^2(\T^d)}
		\leq \sqrt{\frac {2^r}{|\Omega^p_{n/(r+1)}|}} \, \prod_{k=1}^r  \frac{1}{4 |q_k| \eta_k}.
	\end{equation}
\end{lemma}

The assumptions of this lemma are a bit artificial since we cannot assume every hyperplane will have an appropriate $q_k$. However, if $\calX$ locally consists of hyperplanes that satisfy the hypotheses of this lemma for each reference point $x\in \calX$, then the conclusions of \cref{thm:hyperball} and \cref{thm:hypercube} carry over to Fourier matrices modulo constants that do not depend on hyperplane distances or $m$. For instance, the previously considered triangle example $\calX = \{(0,0),(0,\delta),(\delta,0)\}$ falls under this scenario since for each of the three possible reference points, the other two are contained in a hyperplane that is orthogonal to a vector with integer coordinates and whose distance to the reference point is proportional to $\delta$, see Figure \ref{fig:triangleex}. Thus, the machinery developed in this paper  allow us to say that for all sufficiently small $\delta$, we have 
$$
\sigma_{\min}(\Phi_{B_m,\calX})\gtrsim_{s,d} \sqrt{|B_m|_*} (m\delta)  
\andspace \sigma_{\min}(\Phi_{Q_m,\calX})\gtrsim_{s,d} \sqrt{|Q_m|_*} (m\delta). 
$$

In view of this discussion, we present an incomplete strategy and discuss its main bottlenecks. The idea is to perturb an arbitrary $\calX$ to a $\calX^*$ that satisfies the assumptions of \cref{lem:neighborset3}. By using this lemma instead of \cref{lem:neighborset2} in the proofs of \cref{thm:hyperball,thm:hypercube}, we can get an estimate for $\sigma_{\min}(\Phi_{\Omega_m^p,\calX^*})$ in the same spirit as \cref{thm:hyperball,thm:hypercube}. The perturbation from $\calX$ to $\calX^*$ can be achieved through \cref{lem:quanerror}, but unfortunately, the resulting perturbation is too large to be compatible with the next step, which is to relate $\sigma_{\min}(\Phi_{\Omega_m^p,\calX})$ and $\sigma_{\min}(\Phi_{\Omega_m^p,\calX^*})$. This is a classic spectral perturbation problem, but there are limited results for this. A recent preprint \cite{asipchuk2024concerning}, see also \cite{yu2023stability}, examines this problem: in summary, the best results are based on multidimensional Kadec-like theorems that only apply to structured perturbations, while universal perturbation inequalities only hold for very small perturbations. We leave the challenge of deriving estimates for $\sigma_{\min}(\Phi_{\Omega_m^p,\calX})$ based on hyperplane geometry as future work.

\section{Proofs of theorems}
\label{sec:proofthm}

\subsection{Proof of \cref{thm:wellsepcube}}
\label{proof:wellsepcube}

We use \cref{thm:barton} with $\beta d$ acting as $\delta$ to obtain functions $\psi$ and $\varphi$ with properties listed in the referenced theorem. We first concentrate on the analysis of $\calF_{Q_m,\calX}$. Let $\psi_m$ and $\varphi_m$ be dilations such that 
$$
\hat{\psi_m}(\xi)=\hat \psi(m^{-1}\xi) \andspace \hat{\varphi_m}(\xi)= \hat \varphi(m^{-1}\xi).
$$
Note that $\varphi_m$ and $\psi_m$ are supported in $Q_{\beta d/m}$ and that 
$$
\hat{\psi_m}(\xi)\leq \mathbbm{1}_{Q_m}(\xi)\leq \hat{\varphi_m}(\xi) \foreachspace \xi\in\R^d.
$$
Then $\psi_m$ and $\varphi_m$ satisfy the required properties in \cref{lem:poisson} and we have $\Delta_\infty(\calX)\geq \frac {\beta d} m$ by assumption. Thus, we see that
\begin{align*}
	\sqrt{\psi_m(0)}
	\leq \sigma_{\min}(\calF_{Q_m,\calX})
	\leq \sigma_{\max}(\calF_{Q_m,\calX})
	\leq \sqrt{\varphi_m(0)}. 
\end{align*}
It remains to control $\psi_m(0)$ and $\varphi_m(0)$. Using inequalities \eqref{eq:psi0barton} and \eqref{eq:psi0barton2} together with the assumption that $\beta\geq \frac 1 {2\log 2}$, we see that 
\begin{align*}
	\psi_m(0)
	&= m^d \psi(0)
	\geq m^d \big(2-e^{1/(2\beta)} \big)|Q_1| 
	=(2-e^{1/(2\beta)}) |Q_m|, \\
	\varphi_m(0)
	&=m^d \varphi(0)
	\leq m^d e^{1/(2\beta)}|Q_1| 
	=e^{1/(2\beta)} |Q_m|. 
\end{align*}
Combining the above inequalities completes the theorem's proof for $\calF_{Q_m,\calX}$. 

Now we focus on $\Phi_{Q_m,\calX}$. We assume without loss of generality that $m\in \N_+$ since the conclusion remains unchanged if $m$ is not an integer in view of the fact that $|Q_{m}|_*=|Q_{\lfloor m\rfloor}|_*$ and $\Phi_{Q_m,\calX}=\Phi_{Q_{\lfloor m\rfloor},\calX}$. For convenience, set $r=m+1-\epsilon$ where $\epsilon\in (0,1)$ will be adjusted later. Define the dilated functions  $\psi_r$ and $\varphi_r$. Since $r\geq m$, we see that $\varphi_r$ and $\psi_r$ are supported in $Q_{\beta d/r}\subset Q_{\beta d/m}$, that $\Delta_\infty(\calX)\geq \frac{\beta d}m\geq \frac{\beta d}r$, and that 
$$
\hat{\psi_r}(\xi)\leq \mathbbm{1}_{Q_r}(\xi)\leq \hat{\varphi_r}(\xi) \foreachspace\xi\in\R^d.
$$
We have verified that the assumptions of \cref{lem:poisson} are satisfied. Since $m\in\N$ by assumption, while $r=m+1-\epsilon\not\in \Z$, by definition, we have $\Phi_{Q_r,\calX}=\Phi_{Q_m,\calX}$. Applying the lemma, we see that
\begin{align*}
	\sqrt{\psi_r(0)}
	\leq \sigma_{\min}(\Phi_{Q_m,\calX})
	\leq \sigma_{\max}(\Phi_{Q_m,\calX})
	\leq \sqrt{\varphi_r(0)}. 
\end{align*}
It remains to control $\psi_r(0)$ and $\varphi_r(0)$. Repeating the same argument,
\begin{align*}
	\psi_r(0)
	\geq (2-e^{1/(2\beta)})2^d (m+1-\epsilon)^d 
	\andspace
	\varphi_r(0)
	\leq e^{1/(2\beta)}2^d (m+1-\epsilon)^d. 
\end{align*}
Picking $\epsilon=1/2$ and using that $|Q_m|_*=(2m+1)^d$, we obtain
\begin{align*}
	\sigma_{\min}(\Phi_{Q_m,\calX})
	&\geq \sqrt{\psi_r(0)}
	\geq \sqrt{(2-e^{1/(2\beta)})2^d (m+1/2)^d}
	= \sqrt{(2-e^{1/(2\beta)}) |Q_m|_*}, \\
	\sigma_{\max}(\Phi_{Q_m,\calX})
	&\leq \sqrt{\varphi_r(0)}
	\leq \sqrt{e^{1/(2\beta)}2^d (m+1/2)^d}
	= \sqrt{e^{1/(2\beta)} |Q_m|_*}.
\end{align*}

\subsection{Proof of \cref{thm:wellsepball}}

\label{proof:wellsepball}

We only need to prove the theorem for $\alpha\in (\frac{j_{d/2-1,1}}\pi,\frac{j_{d/2,1}}\pi)$ because $c(\alpha)$ is extended to be constant for $\alpha\geq \frac{j_{d/2,1}}\pi$ and the assumption that 
$$
\Delta_2(\calX) \geq \frac{\alpha}m \quad\text{implies} \quad \Delta_2(\calX) \geq \frac{\alpha'}m \forallspace \alpha'\leq \alpha.
$$
Let $\psi$ be the function provided in \cref{thm:goncalves} and define its dilation 
$$
\psi_m(x):=m^d\psi(mx).
$$
Note that $\psi_m$ is supported in $B_{\alpha/m}$ and 
$$
\hat{\psi_m}(\xi)\leq \mathbbm{1}_{B_m}(\xi) \forallspace \xi\in\R^d.
$$
Using \cref{lem:poisson}, we see that 
$$
\sigma_{\min}(\calF_{B_m,\calX})
\geq \sqrt{\psi_m(0)}
= \sqrt{\psi(0) m^d}
= \sqrt{c(\alpha)m^d}. 
$$
We have the same inequality for $\sigma_{\min}(\Phi_{B_m,\calX})$. 

\subsection{Proof of \cref{thm:srball}}
\label{proof:srball}

We fix a $k\in \{1,\dots,s\}$ for now, and consider the decomposition of $\calX$ into the neighborhood set $\calN_2(x_k,\tau,\calX)$ and its complement. We will construct appropriate trigonometric polynomials $b_k$ and $g_k$ such that $b_k(x_k)=g_k(x_k)=1$, $b_k$ vanishes on  $\calN_2(x_k,\tau,\calX)\setminus \{x_k\}$, while $g_k$ vanishes on the complement. 

To construct $g_k$, we will use \cref{prop:localization}, but we first need to check its assumptions hold. By \cref{thm:wellsepball}, for any $\calV\subset \T^d$ such that $\Delta_2(\calV)\geq \tau$, we have
$$
\sigma_{\min}(\Phi_{B_{\alpha/\tau},\calV}) 
\geq \sqrt{c(\alpha) \( \frac{\alpha} \tau \)^d}
= \sqrt{\frac{c(\alpha)}{|B_1|}\frac{|B_{\alpha/\tau}|}{|B_{\alpha/\tau}|_*}} \sqrt{|B_{\alpha/\tau}|_*}.
$$
This together with the assumption $\frac {\alpha\nu} \tau\leq \frac m2$ verifies that the conditions in \cref{prop:localization} hold with $C=\alpha$ and $c$ as the factor in front of $\sqrt{|B_{\alpha/\tau}|_*}$. From the lemma, we obtain a $g_k\in \calP(B_{m/2})$ with the desired interpolation properties and
$$
\|g_k\|_{L^\infty(\T^d)} 
\leq \(\frac{|B_1|}{c(\alpha)}\)^{\nu/2} \(\frac{|B_{\alpha/\tau}|_*}{|B_{\alpha/\tau}|}\)^{\nu/2}. 
$$

For $b_k$, we will use \cref{lem:neighborset1} where the neighborhood set $\calU=\calN_2(x_k,\tau,\calX)-x_k$, $n=\frac m 2$, $r=\nu$, and $p=2$. By the theorem's assumption that $\tau \leq \frac 1{4\sqrt d}$, for any $x_j\in \calN_2(x_k,\tau,\calX)$, we have
$$
|x_j-x_k|_2 \leq \tau \leq \frac 1 {4\sqrt d}.
$$
Also using the assumption $m\geq 4s \sqrt d$, we see that 
$$
\frac m 2 \geq 2 s \sqrt d\geq 2 \nu \sqrt d.
$$
This shows that the assumptions in \cref{lem:neighborset1} are satisfied, from which we conclude (after undoing the translation by $x_k$) that there is a $b_k\in \calP(Q_{m(\nu-1)/(2\nu)})$ such that $b_k(x_k)=1$, $b_k$ vanishes on $\calN_2(x_k,\tau,\calX)\setminus \{x_k\}$ and we have the pointwise estimate 
$$
\|b_k\|_{L^\infty(\T^d)} 
\leq \sqrt{2^{\nu-1}} \prod_{0<|x_j-x_k|_2 \leq \frac \nu m} \, \frac {\nu} {m |x_j-x_k|_2}.
$$

We are ready to conclude. For $\calF_{Q_m,\calX}$, we extend $b_k$ and $g_k$ to elements in $L^\infty(\R^d)$. Let $h$ such that $\hat h= |B_{m/(2\nu)}|^{-1}\mathbbm{1}_{B_{m/(2\nu)}}$. Thus, 
$$
h(0)=1, \quad h\in \calP(B_{m/(2\nu)}),\andspace
\|h\|_{L^2(\R^d)}
=\frac 1 {\sqrt{|B_{m/(2\nu)}|}}. 
$$ 
For $\Phi_{Q_m,\calX}$, we use a different $h$. Consider the Dirichlet kernel, 
$$
h(x) = \frac 1 {|B_{m/(2\nu)}|_*}\sum_{\omega \in B_{m/(2\nu)}\cap \Z^d} e^{2\pi i \omega \cdot x} \andspace
\|h\|_{L^2(\T^d)}
= \frac 1 {\sqrt{|B_{m/(2\nu)}|_*}}.
$$
Clearly $h\in \calP(B_{m/(2\nu)})$ and $h(0)=1$.

For both operators, and for each $k\in \{1,\dots,s\}$, we let $f_k:=h(\cdot -x_k) b_k g_k$. From properties of Fourier transforms, we see that 
$$
\supp(\hat{f_k})
\subset B_{m/(2\nu)} + B_{m(\nu-1)/(2\nu)} + B_{m/2}
\subset B_m.
$$
This shows that $f_k\in \calW(B_m)$ in the case of $\calF_{B_m,\calX}$ and $f_k\in \calP(B_m)$ in the case of $\Phi_{B_m,\calX}$. We also see that 
\begin{align*}
	\|f_k\|_{L^2}
	&\leq \|h\|_{L^2}\|b_k\|_{L^\infty} \|g_k\|_{L^\infty} \\
	&\leq \|h\|_{L^2} \sqrt{2^{\nu-1}}\(\frac{|B_1|}{c(\alpha)}\)^{\nu/2} \(\frac{|B_{\alpha/\tau}|_*}{|B_{\alpha/\tau}|}\)^{\nu/2}  \prod_{0<|x_j-x_k|_2 \leq \frac \nu m} \, \frac {\nu} {m |x_j-x_k|_2}. 
\end{align*}
Applying \cref{lem:duality2} completes the proof of this theorem. 

\subsection{Proof of \cref{thm:srcube}}
\label{proof:srcube}

We fix a $k\in \{1,\dots,s\}$ for now, and consider the decomposition of $\calX$ into the neighborhood set $\calN_\infty(x_k,\tau,\calX)$ and its complement. We will construct appropriate trigonometric polynomials $b_k$ and $g_k$ such that $b_k(x_k)=g_k(x_k)=1$, $b_k$ vanishes on  $\calN_\infty(x_k,\tau,\calX)\setminus \{x_k\}$, while $g_k$ vanishes on the complement. 

To construct $g_k$, we will use \cref{prop:localization}, but we first need to check its assumptions hold. By \cref{thm:wellsepcube}, for any $\calV\subset \T^d$ such that $\Delta_\infty(\calV)\geq \tau$, we have
$$
\sigma_{\min}(\Phi_{Q_{\beta d/\tau},\calV}) \geq \sqrt{(2-e^{1/(2\beta)})|Q_{\beta d/\tau}|_*}.
$$
This together with the assumption $\frac {\beta d \nu} \tau\leq \frac m2$ verifies that the conditions in \cref{prop:localization} hold with $C=\beta d$ and $c= \sqrt{2-e^{1/(2\beta)}}$. From the lemma, we obtain a $g_k\in \calP(Q_{m/2})$ with the desired interpolation properties and
$$
\|g_k\|_{L^\infty(\T^d)} \leq \frac{1}{(2-e^{1/(2\beta)})^{\nu/2}}. 
$$

For $b_k$, we will use \cref{lem:neighborset1} where the neighborhood set $\calU=\calN_\infty(x_k,\tau,\calX)-x_k$, $n=\frac m 2$, $r=\nu$, and $p=\infty$. Using the assumption $\tau \leq \frac 1{4d}$, for any $x_j\in \calN_\infty(x_k,\tau,\calX)$, we have
$$
|x_j-x_k|_1 \leq d \, |x_j-x_k|_\infty \leq d\tau \leq \frac 1 4.
$$
Also using that $m\geq 4s$, we see that 
$
\frac m 2 \geq 2 s\geq 2 \nu.
$
This shows that the assumptions in \cref{lem:neighborset1} are satisfied, from which we conclude (after undoing the translation by $x_k$) that there is a $b_k\in \calP(Q_{m(\nu-1)/(2\nu)})$ such that $b_k(x_k)=1$, $b_k$ vanishes on $\calN_\infty(x_k,\tau,\calX)\setminus \{x_k\}$ and we have the pointwise estimate 
$$
\|b_k\|_{L^\infty(\T^d)} 
\leq \sqrt{2^{\nu-1}} \prod_{0<|x_j-x_k|_1 \leq \frac \nu m} \, \frac {\nu} {m |x_j-x_k|_1}.
$$

We are ready to conclude. For $\calF_{Q_m,\calX}$, we extend $b_k$ and $g_k$ to elements in $L^\infty(\R^d)$. Let $h$ such that $\hat h= |Q_{m/(2\nu)}|^{-1}\mathbbm{1}_{Q_{m/(2\nu)}}$. Thus, 
$$
h(0)=1, \quad h\in \calP(Q_{m/(2\nu)}),\andspace
\|h\|_{L^2(\R^d)}
=\frac 1 {\sqrt{|Q_{m/(2\nu)}|}}. 
$$  
For $\Phi_{Q_m,\calX}$, we use a different $h$. Consider the Dirichlet kernel, 
$$
h(x) = \frac 1 {|Q_{m/(2\nu)}|_*}\sum_{\omega \in Q_{m/(2\nu)}} e^{2\pi i \omega \cdot x} \andspace
\|h\|_{L^2(\T^d)}
= \frac 1 {\sqrt{|Q_{m/(2\nu)}|_*}}.
$$
Clearly $h\in \calP(Q_{m/(2\nu)})$ and $h(0)=1$.

For both operators, and for each $k\in \{1,\dots,s\}$, we let $f_k:=h(\cdot -x_k) b_k g_k$. From properties of Fourier transforms, we see that 
$$
\supp(\hat{f_k})
\subset Q_{m/(2\nu)} + Q_{m(\nu-1)/(2\nu)} + Q_{m/2}
\subset Q_m.
$$
This shows that $f_k\in \calW(Q_m)$ in the case of $\calF_{Q_m,\calX}$ and $f_k\in \calP(Q_m)$ in the case of $\Phi_{Q_m,\calX}$. We also see that 
$$
\|f_k\|_{L^2}
\leq \|h\|_{L^2}\|b_k\|_{L^\infty} \|g_k\|_{L^\infty}
\leq \|h\|_{L^2} \frac{\sqrt{2^{\nu-1}} }{(2-e^{1/(2\beta)})^{\nu/2}}\prod_{0<|x_j-x_k|_1 \leq \frac \nu m} \, \frac {\nu} {m |x_j-x_k|_1}. 
$$
Applying \cref{lem:duality2} completes the proof of this theorem.

\subsection{Proof of \cref{thm:hyperball} and \cref{thm:hypercube}}
\label{proof:hyper}

Both theorems have the same structure, so we prove them simultaneously. Let $p\in \{2,\infty\}$. The starting point of this proof is similar to that of \cref{thm:srball} and \cref{thm:srcube}. For each $k\in \{1,\dots,s\}$, we decompose $\calX$ into a neighborhood set $\calN_p(x_k,\tau,\calX)$ and its complement. Let $g_k\in \calP(\Omega^p_{m/2})$ be the function constructed in \cref{proof:srball} for $p=2$ and \cref{proof:srball} for $p=\infty$. Recall the previously obtained estimates for $\|g_k\|_{L^\infty(\T^d)}$. Now we extend $g_k$ periodically to an element of $L^\infty(\R^d)$.

By assumption, the neighborhood set $\calN_p(x_k,\tau,\calX)$  consists of at most $r$ hyperplanes relative to $x_k$ whose distances to $x_k$ are at least $\eta$. Applying \cref{lem:neighborset2} with $n=\frac m2$, which is justified by the assumption that $\eta\leq \frac{r+1}{2m}$, we obtain a $b_k\in \calW(\Omega^p_{m/2})$ such that $b_k(x_k)=1$, $b_k$ vanishes on $\calN_p(x_k,\tau,\calX)\setminus \{x_k\}$ and 
$$
\|b_k\|_{L^2(\R^d)}
\leq \sqrt{\frac {2^r}{|\Omega^p_{m/(2r+2)}|}} \, \(\frac{r+1}{2m\eta}\)^r.
$$

Define $f_k:=b_kg_k\in L^2(\R^d)$ for each index $k$, we see that $\{f_k\}_{k=1}^s\subset \calW(\Omega_m^p)$ is a family of Lagrange interpolants for $\calX$. Using \cref{lem:duality2} completes this theorem's proof. 

\section{Proofs of propositions and lemmas}
\label{sec:prooflem}

\subsection{Proof of \cref{lem:poisson}}

\label{proof:poisson}

Fix any nonzero $u\in\C^s$. We start with $\calF_{\Omega,\calX}$. Since $\hat\psi$ is integrable, we have 
\begin{align*}
	\|\calF_{\Omega,\calX}u \|_{L^2(\Omega)}^2
	&=\int_{\Omega} \Big| \sum_{k=1}^s u_k e^{-2\pi i\omega \cdot x_k}\Big|^2 \, d\omega 
	\geq \int_{\Omega} \hat\psi(\omega) \Big| \sum_{k=1}^s u_k e^{-2\pi i\omega \cdot x_k}\Big|^2 \, d\omega \\
	&=\sum_{j=1}^s \sum_{k=1}^s \int_{\Omega} \hat\psi(\omega)  \overline{u_j} u_k e^{2\pi i\omega \cdot (x_j-x_k)} \, d\omega 
	= \sum_{j=1}^s \sum_{k=1}^s \overline{u_j} u_k \psi(x_j-x_k). 
\end{align*}
Using that $\psi$ is compactly supported in $B_\delta^p$ and that $\Delta_p(\calX)>\delta$, we see that $\psi(x_j-x_k)=0$ whenever $j\not= k$. This shows that for any $u\in \C^s$, we have 
$$
\|\calF_{\Omega,\calX}u \|_{L^2(\Omega)}^2 \geq \psi(0) |u|_2^2.
$$
Since $u\in \C^s$ is arbitrary, this shows that $\sigma_{\min}(\calF_{\Omega,\calX})\geq \sqrt{\psi(0)}.$ The proof of $\sigma_{\max}(\calF_{\Omega,\calX})\leq \sqrt{\varphi(0)}$ is analogous.

Now we move onto $\Phi_{\Omega,\calX}$. A formal calculation with the Poisson summation formula yields
\begin{equation}
	\label{eq:Phihelp}
	\begin{split}
	&|\Phi_{\Omega,\calX}u|_2^2
	=\sum_{\omega\in\Omega\cap\Z^d} \Big| \sum_{k=1}^s u_k e^{-2\pi i\omega \cdot x_k}\Big|^2
	\geq \sum_{\omega\in\Z^d} \hat\psi(\omega) \Big| \sum_{k=1}^s u_k e^{-2\pi i\omega \cdot x_k}\Big|^2 \\
	&= \sum_{j=1}^s\sum_{k=1}^s \overline{u_j} u_k \sum_{\omega\in\Z^d} \hat\psi(\omega) e^{2\pi i \omega\cdot(x_j-x_k)} 
	=\sum_{j=1}^s\sum_{k=1}^s \overline{u_j} u_k \sum_{\omega\in\Z^d} \psi(\omega+x_j-x_k)
	=\psi(0) |u|_2^2. 
	\end{split}
\end{equation}
We will provide a rigorous derivation at the end, which requires justifying the switch of summation and Poisson summation formula.
Thus, for any unit norm $u\in\C^s$, we have 
$$
|\Phi_{\Omega,\calX}u|_2^2 \geq \psi(0) |u|_2^2. 
$$
This shows that $\sigma_{\min}(\Phi_{\Omega,\calX})\geq \sqrt{\psi(0)}$. For the upper bound on $\sigma_{\max}(\Phi_{\Omega,\calX})$, the argument is analogous, except we use $\varphi$ instead and make appropriate adjustments.

Let us return back to the justification of \eqref{eq:Phihelp}, which are based on standard ideas such as in \cite[Chapter VI, Section 1.15]{katznelson2004introduction}. Let $\eta\colon \R^d\to \R$ be an infinitely differentiable function that is compactly supported in $Q_{1/4}$, $0\leq \hat \eta\leq 1$, and $\hat\eta=1$ on $\Omega$. For any $\epsilon\in (0,1)$, consider the function $\eta_\epsilon(x):=\epsilon^{-d}\eta(x/\epsilon)$. Then $\hat \eta_\epsilon= \eta(\epsilon \omega)$ and in particular, $\hat \eta_\epsilon=1$ on $\Omega$ for all $\epsilon\in (0,1)$. Note $\hat\psi$ is bounded since $\psi$ is continuous and supported in a compact set, and $\hat\eta$ decays faster than any polynomial due to $\eta$ being infinitely differentiable. Hence $\{ \hat\psi(\omega) \hat \eta(\epsilon\omega)\}_{\omega\in\Z^d}\in\ell^1(\Z^d)$, which enables us to switch the order of summation in the following calculation,
\begin{equation}
	\label{eq:Phihelp2}
	\begin{split}
		|\Phi_{\Omega,\calX}u|_2^2
		&=\sum_{\omega\in\Omega\cap\Z^d} \Big| \sum_{k=1}^s u_k e^{-2\pi i\omega \cdot x_k}\Big|^2 \\
		&\geq \sum_{\omega\in\Z^d} \hat\psi(\omega) \hat \eta(\epsilon\omega)\Big| \sum_{k=1}^s u_k e^{-2\pi i\omega \cdot x_k}\Big|^2 \\
		&=\sum_{j=1}^s\sum_{k=1}^s \overline{u_j} u_k \sum_{\omega\in\Z^d} \hat\psi(\omega) \hat \eta(\epsilon\omega) e^{2\pi i \omega\cdot(x_j-x_k)}.
	\end{split}
\end{equation}
Due the assumption that $\psi$ is supported in $\Omega_\delta^p$ for some $\delta\in (0,\frac 12)$, we see $\psi$ can be interpreted as a continuous function on $\T^d$ and that $\hat\psi(\omega)$ is the Fourier coefficient of $\psi$ at $\omega\in \Z^d$. Likewise, $\eta_\epsilon$ is compactly supported in $Q_{\epsilon/4}$, so it is also an infinitely differentiable function on $\T^d$ and $\hat \eta(\epsilon\omega)$ is the Fourier coefficient of $\eta_\epsilon$ at $\omega\in \Z^d$. Thus, the sum over $\omega\in \Z^d$ on the right hand side of \eqref{eq:Phihelp2} is the Fourier series of $\psi*\eta_\epsilon$ (treated as a periodic function) evaluated at $x_j-x_k$. This Fourier series converges pointwise due to $\{ \hat\psi(\omega) \hat \eta(\epsilon\omega)\}_{\omega\in\Z^d}\in\ell^1(\Z^d)$, and so
$$
\sum_{j=1}^s\sum_{k=1}^s \overline{u_j} u_k\sum_{\omega\in\Z^d} \hat\psi(\omega) \hat \eta(\epsilon\omega) e^{2\pi i \omega\cdot(x_j-x_k)} 
= \sum_{j=1}^s\sum_{k=1}^s \overline{u_j} u_k (\psi*\eta_\epsilon)(x_j-x_k).
$$
Combining the above observations yields the inequality 
$$
|\Phi_{\Omega,\calX}u|_2^2
\geq \sum_{j=1}^s\sum_{k=1}^s \overline{u_j} u_k (\psi*\eta_\epsilon)(x_j-x_k). 
$$
Now we are in position to take the limit as $\epsilon\to 0$. By assumption, $\eta$ is integrable and $\int_{\T^d} \eta_\epsilon=\int_{\T^d} \eta=\hat\eta(0)=1$. Then $\{\eta_\epsilon\}_{\epsilon>0}$ is an approximation to the identity as $\epsilon\to 0$. Hence, $\psi*\eta_\epsilon\to \psi$ pointwise since $\psi$ is assumed to be continuous. Therefore, taking the limit now shows that 
$$
|\Phi_{\Omega,\calX}u|_2^2
\geq \sum_{j=1}^s\sum_{k=1}^s \overline{u_j} u_k \psi(x_j-x_k)
=\sum_{j=1}^s |u_j|^2 \psi(0)
=\psi(0)|u|^2_2. 
$$
This now proves inequality \eqref{eq:Phihelp}.

\subsection{Proof of \cref{prop:duality}}
\label{proof:duality}

The proof of the two statements are similar, so we only prove the first one for $\calF_{\Omega,\calX}$. Let $\mu = \sum_{k=1}^s v_k \delta_{x_k}$. On one hand, using the interpolation properties of $f$, we have 
$$
\Big|\int_{\T^d} \overline f \, d\mu\Big| 
= \Big| \sum_{k=1}^s v_k \overline{f(x_k)} \Big|  
=  \Big| |v|_2^2 + \sum_{k=1}^s v_k \overline{\epsilon_k} \Big| 
\geq |v|_2^2 - |v|_2 |\epsilon|_2
= 1-|\epsilon|_2. 
$$
On the other hand, using that $f\in \calW(\Omega)$ and Cauchy-Schwarz, we have
\begin{align*}
	\Big|\int_{\T^d} \overline f(x) \, d\mu(x) \Big|
	&= \Big| \int_{\T^d} \int_{\Omega} \overline {\hat f(\omega)}e^{-2\pi i \omega\cdot x} \, d\omega d\mu(x) \Big| \\
	&= \Big|\int_{\Omega} \overline {\hat f(\omega)} \hat\mu(\omega) \, d\omega \Big|
	\leq \|\hat f\|_{L^2(\Omega)} \|\hat\mu\|_{L^2(\Omega)}
	= \|f\|_{L^2(\R^d)} \|\calF_{\Omega,\calX} v\|_{L^2(\Omega)}.
\end{align*}
Combining the above and rearranging completes the proof. 

\subsection{Proof of \cref{lem:duality2}}
\label{proof:duality2}

If $\{f_k\}_{k=1}^s$ is a family of Lagrange interpolants for $\calX$, then $f=\sum_{k=1}^s v_kf_k$ interpolates a unit $\ell^2$ norm $v$ on $\calX$ and 
$$
\|f\|_{L^2}
\leq \sum_{k=1}^s |v_k| \|f_k\|_{L^2}
\leq |v|_1 \max_{1\leq k\leq s} \|f_k\|_{L^2} 
\leq \sqrt s \max_{1\leq k\leq s} \|f_k\|_{L^2}.
$$ 
Applying \cref{prop:duality} where $v$ is a singular vector corresponding ot the smallest singular value completes the proof.

\subsection{Proof of \cref{prop:localization}}
\label{proof:localization} 

We fix a $k\in \{1,\dots,s\}$ for now and let $\calG_k:=\calX\setminus \calN_p(x_k,\tau,\calX)$. We apply \cref{prop:decomp} to $\calG_k$ to obtain a disjoint union
$$
\calX = \calN_p(x_k,\tau,\calX) \cup \bigcup_{\ell=1}^{\nu_p(\tau,\calG_k)}\calG_{k,\ell} \wherespace \Delta_p(\calG_{k,\ell})>\tau \foreachspace \ell \in \{1,\dots,\nu_p(\tau,\calG_k)\}.
$$

By definition of the neighborhood set, we have $\Delta_p(\calG_{k,\ell}\cup \{x_k\})>\tau$ for each $\ell$ as well. Employing the assumed inequality \eqref{eq:wellsepbound} for $\calG_{k,\ell}\cup \{x_k\}$ acting as $\calW$, we obtain 
$$
\sigma_{\min}\big(\Phi_{\Omega^p_{C/\tau},\, \calG_{k,\ell}\cup \{x_k\}}\big)
\geq c \sqrt{|\Omega^p_{C/\tau}|_*} \foreachspace \ell\in \{1,\dots,\nu_p(\tau,\calG_k)\}. 
$$

For each $\ell$, we apply \cref{prop:interpolation2} to $\calG_{k,\ell} \cup \{x_k\}$ to obtain the existence of a polynomial $g_{k,\ell} \in \calP(\Omega^p_{C/\tau})$ such that $g_{k,\ell}(x_k)=1$, $g_{k,\ell}$ vanishes on $\calG_{k,\ell}$, and  
$$
\|g_{k,\ell}\|_{L^\infty(\T^d)}
\leq \frac {\sqrt{|\Omega^p_{C/\tau}|_*}} {\sigma_{\min}\big(\Phi_{\Omega^p_{C/\tau}, \calG_{k,\ell} \cup \{x_k\}}\big)}
\leq \frac 1 {c}. 
$$
Doing this for each $\ell$ and multiplying these $\nu_p(\tau,\calG_k)$ polynomials together provides us with the desired $g_k$. There are $\nu_p(\tau,\calG_k)$ many $g_{k,\ell}$, and note that $\calG_k\subset\calX$ implies $\nu_p(\tau,\calG_k)\leq \nu_p(\tau,\calX)$. Thus, we see that 
$$
\|g_k\|_{L^\infty(\T^d)}
\leq \frac 1 {c^{\nu_p(\tau,\calG_k)}}
\leq \frac 1 {c^{\nu_p(\tau,\calX)}}.
$$
As for the support of $\hat{g_k}$, note that each $\hat{g_{k,\ell}}$ is supported in $\Omega^p_{C/\tau}\cap \Z^d$. There are at $\nu_p(\tau,\calG_k)$ many functions and recall the sparsity assumption \eqref{eq:densitycondition}. Thus, we have
$$
\supp(\hat{g_k})
\subset \Omega^p_{C/\tau} + \cdots + \Omega^p_{C/\tau} 
\subset \Omega^p_{C\nu_p(\tau,\calG_k) /\tau}
\subset \Omega^p_{C\nu_p(\tau,\calX) /\tau}
\subset \Omega^p_{m/2}. 
$$
This completes the first assertion's proof about existence of $g_k$.

For the second assertion, by \cref{prop:interpolation2} and the assumption that $\sigma_{\min}(\Phi_{\Omega_{m/2}^p,\calN_p(x_k,\tau,\calX)})>0$, we can obtain the existence of a $f_k\in \calP(\Omega_{m/2}^p)$ such that $f_k(x_k)=1$, $f_k$ vanishes on $\calN_p(x_k,\tau,\calX)\setminus \{x_k\}$, and 
$$
\|f_k\|_{L^2(\T^d)}
\leq \frac {1}{\sigma_{\min}(\Phi_{\Omega_{m/2}^p,\calN_p(x_k,\tau,\calX)})}.
$$ 
From here, we see that $\{f_kg_k\}_{k=1}^s\subset \calP(\Omega_m^p)$ is a family of Lagrange polynomials for $\calX$. The proof is complete once we use the above estimates for $\|g_k\|_{L^\infty(\T^d)}$ and $\|f_k\|_{L^2(\T^d)}$ together with \cref{lem:duality2}.

\subsection{Proof of \cref{lem:quanerror}}
\label{proof:quanerror} 

For $a\in \R$, we denote $[a]\in \Z$ by the unique integer such that $|a-[a]|<1$ and $|[a]|\leq |a|$. We extend this operation to vectors in the following way. For $a\in\R^d$, we let $[a]\in\Z^d$ such that $[a]_k=[a_k]$ for each $k\in \{1,\dots,d\}$. We have the trivial inequalities: for all $a\in \R^d$ and $p\in [1,\infty]$, we have
\begin{equation}
	\label{eq:quanerror0}
	|a-[a]|_\infty< 1 \andspace |[a]|_p\leq |a|_p. 
\end{equation}

Let $v\in \R^d$ such that $|v|_p=1$ and $|v\cdot u|=|u|_{p'}$. That is, $v$ is the unit $\ell^p$ dual of $u\in \ell^{p'}$ and explicitly, $v_k={|u_k|^{p'-1} \sign(u_k)}/{|u|_{p'}^{p'-1}}$ with the convention that $\sign(u_k)=0$ if $u_k=0$. From here, we set 
$$
a:= \frac{v}{2\alpha} \andspace q:= [a].
$$ 
The first inequality \eqref{eq:lemquanerror} of this lemma follows immediately from \eqref{eq:quanerror0} because 
$
|q|_p=|[a]|_p
\leq |a|_p =\frac 1 {2\alpha}.
$
Moving on, we first use H\"older's inequality, we get
$$
|q\cdot u|\leq |q|_p |u|_{p'}\leq \frac 1 {2\alpha} |u|_{p'}\leq \frac 1 2. 
$$
To prove the claimed lower bound for $|q\cdot u|$, we use \eqref{eq:quanerror0} to see that 
$
|(q-a)\cdot u|
\leq |q-a|_\infty |u|_1
\leq |u|_1
\leq d^{1/p} |u|_{p'}. 
$
Using this inequality, that $a\cdot u = \frac {1} {2\alpha}|u|_{p'}$, and the assumption $\alpha \leq \frac {1}{4 d^{1/p}}$, we see that
$$
|q\cdot u|
\geq |a\cdot u| - |(q-a)\cdot u|
\geq \frac 1 {2\alpha}|u|_{p'} - d^{1/p} |u|_{p'}
\geq \frac 1 {4\alpha} |u|_{p'}. 
$$
This proves the middle two inequalities in \eqref{eq:lemquanerror} of this lemma. 

Now we move onto the final inequality. We define the sinc kernel, $\sinc(t):=\frac{ \sin (\pi t)}{\pi t}$. This function will naturally appear in our analysis because for all $|t|\leq \frac 1 2$, we have
\begin{equation}
	\label{eq:lawcos}
	|1-e^{2\pi it}|
	=\sqrt{2-2\cos(2\pi t)}
	= 2 \pi t \, \sinc(t).
\end{equation}
Using the lower bound for $|q\cdot u|$ in \eqref{eq:lemquanerror} and \eqref{eq:lawcos}, we see that 
$$
|1-e^{2\pi i q\cdot u}| 
= 2 \pi |q\cdot u| \, \sinc(|q\cdot u|)
\geq \frac{\pi}{2 \alpha} |u|_{p'} \, \sinc\(\frac{|u|_{p'}}{4\alpha}\)
\geq \frac{\pi}{2 \alpha} |u|_{p'} \, \sinc\(\frac 1 4\)
= \frac{\sqrt 2}{\alpha} |u|_{p'}, 
$$
where for the final step, we used that $\sinc$ is decreasing away from zero on $[-\frac 12,\frac 12]$ and that $\frac{|u|_{p'}}{4\alpha}\leq \frac 1 4$. This completes the proof of the lemma. 

\subsection{Proof of \cref{lem:neighborset1}}

\label{proof:neighborset1}

For $\calU=\{0\}$, there is nothing to prove since we can set $f=1$ and it satisfies the claimed properties and inequality \eqref{eq:neighbor1}, where we recall that a product over an empty set is defined to be 1. From here onward, assume that $|\calU|\geq 2$. Define the subsets
$$
\calI:=\Big\{u \in \calU \colon 0<|u|_{p'} \leq \frac {r} {2n} \Big\} \andspace \calJ:=\calU\setminus(\calI\cup\{0\}). 
$$

We will first construct a function $h$ such that $h(0)=1$ and $h$ vanishes on $\calI$. If $\calI=\emptyset$, then we set $h=1$. Now suppose that $\calI\not=\emptyset$. For each $u\in\calI$, we use \cref{lem:quanerror} where $\frac r {2n}$ plays the role of $\alpha$. This is justified since the lemma assumes that $|u|_{p'}\leq \frac r {2n}$ and $n\geq 2d^{1/p}r$. Doing so, for each $u\in \calI$, there is a $q(u)\in\Z^d$ such that  
\begin{equation}
	\label{eq:quanerror}
	|q(u)|_p\leq \frac{n}{r}, \quad 
	\frac r {2n} |u|_{p'} \leq |q(u)\cdot u| \leq \frac 1 2, \andspace 
	|1-e^{2\pi i q(u) \cdot u}| \geq \frac{2 \sqrt 2 \, n}{r} |u|_{p'}.  
\end{equation}
We define $h$ by the formula,
$$
h(x)
:=\prod_{u\in\calI} \frac{e^{2\pi i q (u)\cdot x}-e^{2\pi iq(u)\cdot u}}{1-e^{2\pi iq(u)\cdot u}}. 
$$
Note that $h$ is well defined since the second inequality in \eqref{eq:quanerror} implies each denominator is nonzero. By construction $h(0)=1$ and $h=0$ on $\calI$. Next, due to the first inequality in \eqref{eq:quanerror}, each function in the product term in the definition of $h$ belongs to $\calP(\Omega_{n/r})$. Consequently,
$$
\supp(\hat h)
\subset \Omega_{n/r}+\cdots + \Omega_{n/r}
\subset \Omega_{|\calI|n/r}.
$$	
To control $\|h\|_{L^\infty(\T^d)}$, we use \eqref{eq:quanerror} to see that
\begin{equation}
	\label{eq:hhelp2}
	\|h\|_{L^\infty(\T^d)}
	\leq \prod_{u\in\calI} \frac{2}{|1-e^{2\pi i q(u) \cdot u}|}
	\leq 2^{|\calI|/2} \prod_{u\in \calI} \, \frac {r} {2n |u|_{p'}}. \\
\end{equation}
Hence, regardless of whether $\calI$ is empty or not, we have found a $h\in \calP(\Omega_{|\calI|/r})$ with the desired interpolation properties, and satisfies estimate \eqref{eq:hhelp2}.

Now we construct a function $g\in\calP(\Omega_{|\calJ|n/r})$ such that $g(0)=1$ and $g$ vanishes on $\calJ$, and $\|g\|_{L^\infty(\T^d)} \leq 2^{|\calJ|/2}.$ If $\calJ=\emptyset$, then we set $g=1$ and we are done. Suppose $\calJ\not=\emptyset$. For each $u\in \calJ$, we apply \cref{lem:quanerror} where $|u|_{p'}$ plays the role of $\alpha$, which is justified by the assumption that $|u|_{p'}\leq \frac 1  {4d^{1/p}}$. Then there is a $q(u)\in \Z^d$ such that 
\begin{equation}
	\label{eq:quanerror2}
	|q(u)|_p\leq \frac{1}{2|u|_{p'}}, \quad 
	\frac 1 {4} \leq |q(u)\cdot u| \leq \frac 1 2, \andspace 
	|1-e^{2\pi i q(u) \cdot u}| \geq \sqrt 2.  
\end{equation}
With these considerations completed, we define the function 
$$
g(x)=\prod_{u\in\calJ} \frac{e^{2\pi i q(u)\cdot x}-e^{2\pi i q(u) \cdot u}}{1-e^{2\pi i q(u) \cdot u}}. 
$$
This is well defined since each term in the denominator is positive in view of the second inequality in \eqref{eq:quanerror2}. By construction, $g(0)=1$ and $g$ vanishes on $\calJ$. Due to the first inequality in \eqref{eq:quanerror2} and definition of $\calJ$, we have $|q(u)|_p< \frac n r$, and consequently, $g\in \calP(\Omega_{|\calJ|n/r})$. Using the final inequality in \eqref{eq:quanerror2}, we obtain the pointwise estimate
\begin{equation}
	\label{eq:ghelp1}
	\|g\|_{L^\infty(\T^d)}
	\leq \prod_{u\in \calJ} \frac{2}{|1-e^{2\pi i q(u) \cdot u}|}
	\leq 2^{|\calJ|/2}.
\end{equation}

To finish the proof, we let $f:=gh$ which enjoys the interpolation properties that $f(0)=1$ and $f=0$ on $\calU\setminus\{0\}$. The estimate for $\|f\|_{L^\infty(\T^d)}$ follows by combining our previous estimates  \eqref{eq:hhelp2} and \eqref{eq:ghelp1}. Finally, we have $f\in \calP(\Omega_{n(r-1)/r})$ because $|\calU|\leq r$ by assumption and $\calU = \{0\}\cup \calI \cup \calJ$ is a disjoint union, so 
$$
\supp(\hat f)
\subset \supp(\hat g) + \supp(\hat h)
\subset \Omega_{|\calI|n/r} + \Omega_{|\calJ|n/r}
\subset \Omega_{(|\calU|-1)n/r}
\subset \Omega_{n(r-1)/r}.
$$	

\subsection{Proof of \cref{lem:neighborset2}}

\label{proof:neighborset2}

Enumerate the $r$ hyperplanes by $\calH_1,\dots,\calH_r$. Since $\eta_k$ is the $\ell^2$ distance between $\calH_k$ and $0$, we can find a unit $\ell^2$ norm vector  $\theta_k\in \R^d$ that is orthogonal to $\calH_k$ and  $u_k:=\eta_k \theta_k\in \calH_k$. Next, we consider the function $h\in L^\infty(\R^d)$ by 
$$
h(x):=\prod_{k=1}^r \frac{e^{2\pi i x\cdot n\theta_k/(r+1)}-e^{2\pi i u_k\cdot n\theta_k/(r+1)}}{1-e^{2\pi i u_k\cdot n\theta_k/(r+1)}}. 
$$
We make several comments about $h$. First, the denominator is always nonzero because 
\begin{equation}
	\label{eq:thetahelp}
	\Big|u_k\cdot \frac n{r+1} \theta_k \Big| = \frac n{r+1}  \eta_k
	\leq \frac 14,
\end{equation}
where the last inequality follows by assumption. Consequently, $h$ is well defined. Second, the $k$-th factor in the definition of $h$ is constant on hyperplanes parallel to $\calH_k$, i.e., orthogonal to $\theta_k$. Indeed, for any $x\in \R^d$ and $v\in \R^d$ that is orthogonal to $\theta_k$, we see that $(x+v)\cdot\theta_k = x\cdot \theta_k$. Since the $k$-th factor is zero when evaluated at $u_k\in \calH_k$, we see that the $k$-th factor is zero on $\calH_k$. Thus, $h(0)=1$ and $h$ vanishes on $\calH_1\cup \cdots \cup \calH_r$. Finally, in the sense of distributions, the Fourier transform of $h$ is supported in the set 
\begin{equation}
	\label{eq:hFouriersupp}
	\left\{ \frac n {r+1} \theta_1, \, \frac n {r+1} \theta_2, \, \, \dots, \, \frac n {r+1} \theta_r \right\}.
\end{equation}

Next, define $g$ such that $\hat g= |\Omega^p_{n/(r+1)}|^{-1} \mathbbm{1}_{\Omega^p_{n/(r+1)}}$. Note that
\begin{equation}
	\label{eq:g2norm}
	g(0)=1, \quad g\in \calW(\Omega^p_{n/(r+1)}) \andspace
	\|g\|_{L^2(\R^d)}= \frac 1 {\sqrt{|\Omega^p_{n/(r+1)}|}}. 
\end{equation}
Finally, we define $f:=gh$, which is an element of $L^2(\R^d)$. In view of \eqref{eq:hFouriersupp}, the Fourier transform of $f$ is contained in
\begin{align*}
	&\Omega^p_{n/(r+1)} +\frac n {r+1} \theta_1 + \frac n {r+1} \theta_2 + \cdots + \frac n {r+1} \theta_r \\
	&\quad \quad\quad \subset \Omega^p_{n/(r+1)} + B_{n/(r+1)} + \cdots + B_{n/(r+1)}
	\subset \Omega^p_n,
\end{align*}
where we used that $\Omega^2_{n/(r+1)}\subset \Omega^p_{n/(r+1)}$ since $p\in [2,\infty]$. This establishes that $f\in \calW(\Omega^p_n)$ and has the claimed interpolation properties. 

Next, we concentrate on bounding $\|h\|_{L^\infty(\R^d)}$. Using \eqref{eq:thetahelp}, we have
$$
|1-e^{2\pi i u_k\cdot n\theta_k/(r+1)}|
=|1-e^{2\pi i n\eta_k/(r+1)}|
\geq \frac{ 2\pi n \eta_k }{r+1}\, \sinc \(\frac 14\)
= \frac{ 4 \sqrt 2 n \eta_k }{r+1}.
$$
From this inequality, we see that 
\begin{equation}
	\label{eq:hnorm}
	\|h\|_{L^\infty(\R^d)}
	\leq \prod_{k=1}^r \frac{2}{|1-e^{2\pi i u_k\cdot n\theta_k/(r+1)}|}
	\leq \prod_{k=1}^r \frac{r+1}{2 \sqrt 2 n\eta_k}
	= 2^{r/2} \prod_{k=1}^r \frac{r+1}{4 n\eta_k}. 
\end{equation}
The proof is completed once we use H\"older's inequality $\|f\|_{L^2(\R^d)}\leq \|g\|_{L^2(\R^d)} \|h\|_{L^\infty(\R^d)}$ with inequalities \eqref{eq:g2norm} and \eqref{eq:hnorm}. 

\subsection{Proof of \cref{lem:neighborset3}}
\label{proof:neighborset3}

Enumerate the $r$ hyperplanes by $\calH_1,\dots,\calH_r$. Since $\eta_k$ is the $\ell^2$ distance between $\calH_k$ and $0$, we can find a unit $\ell^2$ norm vector  $\theta_k\in \R^d$ that is orthogonal to $\calH_k$ and  $u_k:=\eta_k \theta_k\in \calH_k$.  Consider the trigonometric polynomial $h$, where
$$
h(x):=\prod_{k=1}^r \frac{e^{2\pi i q_k \cdot x}-e^{2\pi i q_k\cdot u_k}}{1-e^{2\pi i q_k\cdot u_k}}. 
$$
The denominators of this function are all nonzero because $\theta_k$ is parallel to $q_k$, which establishes that $|q_k \cdot u_k| = |q_k|_2 |\theta_k\cdot u_k| = \eta_k |q_k|_2 >0$. To control the denominators, notice that $\eta_k |q_k|_2\leq \eta_k \frac n {r+1}\leq \frac 14$, so 
$$
|1-e^{2\pi i q_k\cdot u_k}|
\geq 2\pi \eta_k |q_k|_2 \sinc\(\frac 14\)
= 4 \sqrt 2 \eta_k |q_k|_2.
$$
Next, let $h$ be the normalized Dirichlet kernel for $\Omega^p_{n/(r+1)}$. Finally, we define $f=gh$ which satisfies the desired interpolation inequalities, and we have 
$$
\|f\|_{L^2(\T^d)}
\leq \|g\|_{L^2(\T^d)} \|h\|_{L^\infty(\T^d)}
\leq \frac 1 {\sqrt{|\Omega^p_{n/(r+1)}|_*}} 2^{r/2} \prod_{k=1}^r \frac{1}{4 |q_k| \eta_k}.
$$
Finally, we see that $f\in \calP(\Omega^p_n)$ because $g\in \calP(\Omega^p_{n/(r+1)})$ and $|q_k|_2\leq \frac n {r+1}$ for each $k\in \{1,\dots,r\}$ implies $h\in \calP(B_{nr/(r+1)})$. 

\section*{Conclusion}

In this paper, we studied the smallest singular values of a continuous bandlimited Fourier operator and a discrete counterpart. We derived estimates that show how the choice of metric and geometry of $\calX$ relative to the bandwidth $2m$ influence their singular values. The three types of estimates proceed in increasing complexity of the geometric assumptions on $\calX$, from minimum separation, to local sparsity, and to hyperplane decompositions.

There are numerous extensions and adaptations of our results. The main theorems of this paper are stated for $B_m$ and $Q_m$. They can be adapted to $\Omega_m^p$ for all $p\in [1,\infty]$, but doing so will introduce additional (and possibly unnatural) dependence on the dimension. For instance, minorants for the ball \cite{gonccalves2018note} and cube \cite{barton1999analogs}, can be adapted to $\Omega^p_m$, but one pays a price in dimension due to H\"older's inequality. We are unaware of precise results for $\ell^p$ ball minorants.  Nonetheless, the lemmas in \cref{sec:sr1} and \cref{sec:sr2} hold for all $p\in [1,\infty]$ and $p\in [2,\infty]$ respectively. 

One could consider sampling the Fourier transform from $TB_m$ or $TQ_m$ for some invertible linear transformation $T$ on $\R^d$. This would correspond to using a (possibly rotated) ellipsoidal or parallelpiped aperture, respectively. In this case, our results for the continuous Fourier operator readily generalize due to simple transform properties -- one would just replace the $\ell^p$ distance on $\calX$ with the $\ell^p$ metric in the transformed space, $|x|_{p,T}:=|Tx|_p$. For Fourier matrices, we expect exactly the same effect. However, it does not have nice transformation properties like its continuous counterpart (e.g., $B_m$ is closed under rotations while $B_m\cap\Z^d$ is not), so we cannot simply perform a change of variables -- one would need to modify some of the lemmas and arguments in this paper.

A natural application of our results is to multidimensional super-resolution, where the problem is to recover $\mu$ from a perturbed Fourier transform restricted to $\Omega\subset \R^d$ or $\Omega\cap\Z^d$. Our results have two main implications. A first application is the performance of multidimensional super-resolution algorithms. The currently available results for multidimensional MUSIC \cite{liao2015music} can be greatly improved using the theorems in this paper. There are significant challenges in analyzing multidimensional extensions of ESPRIT and MPM. 

A second application is to understand optimal rates of recovery. Essentially, if there are two measures, supported in $\calX$ and $\calX^*$ and with amplitudes $u$ and $u^*$ respectively, then the distance between their measurements is 
$$
\|\calF_{\Omega^p_m,\calX} u - \calF_{\Omega^p_m,\calX^*} u^* \|_{L^2(\Omega)}
=\|\calF_{\Omega^p_m,\calX\cup\calX^*}(u,u^*)\|_{L^2(\Omega)}. 
$$
These measures are indistinguishable from data if this term is less than the noise level measured in $L^2(\Omega)$, and this term is precisely the smallest singular value of $\calF_{\Omega^p_m,\calX\cup\calX^*}$ when $(u,u^*)$ is a corresponding right singular vector. This provides the connection between optimal recovery rates for super-resolution and the smallest singular value of Fourier operators, which has been previously exploited \cite{donoho1992superresolution,demanet2015recoverability,li2021stable,batenkov2021super,liu2021mathematical}.

It is worth discussing implications of the results in \cref{sec:sr1} versus \cref{sec:sr2}. When do we expect to encounter worst case or generic sets? For certain imaging application such as astronomy, we generally do not expect stars to be densely arranged on lines or curves, so the hyperplane results are more appropriate. This is a blessing of dimensionality in the sense that super-resolution is more robust to noise than what is suggested by worst case analysis or one-dimensional results. On the other hand, for other applications like 2D super-resolution florescence microscopy, the chemical compounds that stochastically emit light attach to the boundary of cells which are locally quadratic curves. If many of them simultaneously activate, then we expect to encounter the worst case situation.

\section*{Acknowledgments}

WL is supported by NSF-DMS Award \#2309602, a PSC-CUNY grant, and a start-up fund from the Foundation for City College. WL thanks Dmitry Batenkov, Jacob Carruth, Jacky Chong, Albert Fannjiang, Felipe Gon\c calves, and Wenjing Liao for valuable and insightful discussions.

\bibliography{multidimensionalFourierBib.bib}
\bibliographystyle{plain}

\end{document}